\documentclass[a4paper]{article}
\usepackage{graphicx} % Required for inserting images
\usepackage{xcolor}
\usepackage{tabularx, booktabs}
\usepackage{subcaption}
\usepackage{multirow}
\usepackage{natbib}
\usepackage{multicol}
\usepackage{amsmath}
\usepackage{amssymb}
\usepackage{todonotes}
\begin{document}
\begin{center}{
\Large{Learning Dominant States in Elementary Resource Constrained Shortest Path Problems}}
\end{center}

%Authors
\begin{center}
\index{Author1}Saverio Basso$^{1}$, \index{Author2}Matteo Salani$^{1}$\\ 
\vspace{2ex}$^{1}$Dalle Molle Institute for Artificial Intelligence (IDSIA), USI-SUPSI, Lugano, Switzerland \\ 
\end{center}
\vspace{3ex}
%\maketitle

%\author*[1]{\fnm{Saverio} \sur{Basso}}\email{saverio.basso@supsi.ch}
%\equalcont{These authors contributed equally to this work.}

%\author[1]{\fnm{Matteo} \sur{Salani}}\email{matteo.salani@supsi.ch}
%\equalcont{These authors contributed equally to this work.}

%\affil[1]{\orgname{SUPSI, Dalle Molle Institute for Artificial Intelligence (IDSIA)}, \orgaddress{\street{Via la Santa 1}, \city{Lugano-Viganello}, \postcode{CH-6962}, \country{Switzerland}}}

\abstract{
	In this work, we investigate whether machine learning can be leveraged to identify promising states in dynamic programming algorithms, focusing on Elementary Resource Constrained Shortest Path Problems (ERCSPP).

    More in detail, we solved 41 single resource instances from SPPRCLIB using iterative relaxation techniques through the PathWyse library, systematically collecting all generated states (i.e. labels). We designed ad-hoc features computable in constant time and constructed two datasets: one containing all generated labels (G) and another with only those inserted into data pools (I), totaling several hundred million labels.

    Machine learning tools are then employed to explore these datasets, revealing significant patterns between successive relaxations. Leveraging these insights, we propose a normalization approach and apply supervised learning techniques to distinguish dominating states, both within subsequent relaxations of the same problem and in previously unseen instances. Our results demonstrate the effectiveness of this approach on Dataset G, while for Dataset I, performance varies, showing strong results within the same instance but declining for unseen ones. Overall, these findings open new perspectives for the development of data-driven dynamic programming algorithms.
}

%\keywords{ESPPRC, machine learning, dynamic programming, dominance}

%\maketitle

\section{Introduction}

The shortest path problem is one of the most studied combinatorial problems. While original attempts can be dated back to the 19th century in methods to solve path finding in mazes  (\cite{lucas1882recreations} and \cite{tarry1895probleme}), modern mathematical study of the shortest path problem began in the mid-20th century with the work by \cite{dijkstra1959note}. The interested reader can discover the history of the shortest path problem in \cite{schrijver2012history}.

The shortest path problem is efficiently solvable in polynomial time. However, the introduction of resource constraints transforms it into a significantly more complex and challenging problem, in the class of NP-Hard problems. These constraints can include limits on time, cost, or other resources that must be managed along the path.
Another layer of complexity is added when the underlying network contains negative cost cycles. These cycles can lead to situations where the total path cost can be reduced by traversing the cycle repeatedly. The solution to this variant asks for the computation of the \textit{Elementary} Shortest Path Problem with Resource Constraints (ESPPRC).
There are several practically relevant applications of the ESPPRC in transportation, telecommunication, finance, and scheduling. Most notably, the ESPPRC is found as subproblem in routing and workforce scheduling problems solved by column generation.

%% Introduzione ML
Recently, data-driven algorithms have emerged as complementary tools to tackle combinatorial optimization problems. In particular, leveraging on large datasets, Machine Learning (ML) approaches are capable of effectively driving search algorithms in the very large solution space of combinatorial problems. While pioneering work in the direction of combining learning and optimization dates back to the early 1990s, see for example \cite{dorigo1992optimization}, this research avenue has been reshaped more recently towards the integration of ML to drive inner decisions of more complex algorithms for combinatorial optimization, \cite{BENGIO2021}. 

In this work, we focus on the ESPPRC and on one of the state of the art algorithms to solve the problem to proven optimality. Namely, we focus on the Decremental State Space Relaxation algorithm (DSSR) \citep{RighiniSalani2008} which is based on a dynamic programming approach capable of effectively explore the large state space associated to the problem. We intend to explore the effectiveness of ML techniques to identify promising states (later on called \textit{Pareto states}) that are worth exploring first.
%extend here?
In fact, the performance of the DSSR algorithm is strongly influenced by the state space exploration policy. For example, Table \ref{tab:policy} illustrates the number of states explored (Attempts) and the percentage of those inserted in data pools (Success) under two different strategies when solving test-bed ESPPRC instances with one resource constraint. The node policy explores all the unexplored states associated with a node of the network before moving to the next node, while the greedy policy explores the most promising state in the whole network first. Overall, the node policy proved almost three times faster than the greedy policy in all instances. Indeed, these findings further motivate the investigation of data-driven techniques to identify Pareto optimal states as they are generated and explored by the DSSR algorithm.

\begin{table}[h!]
\small

\begin{tabularx}{\textwidth}{X*{6}{r}}\toprule
%        &           &         &&  &\multicolumn{2}{c}{Insert} \\
%\cmidrule(lr){6-7} 
Policy	&Time [s]	&Fastest	&	&Attempts	&Success\\
	\midrule

Node policy	&50.52	&42	&	&1697524652	&2.89\%\\
	&	&	&	&	&	&\\
Greedy policy	&146.16	&0	&	&5521806107	&1.85\%\\
\addlinespace
\bottomrule
\end{tabularx}
\caption[Extension policy profiling]
{Labeling algorithm profiling over different extension policies when solving ESPPRC instances.} \label{tab:policy}
\end{table}

The paper is organized as follows: in Section \ref{sec:literature} we briefly review the most relevant contributions to solve the (E)SPPRC. Section \ref{sec:algo} recalls the details of the selected labeling algorithm for the problem and defines fundamental research questions. The generation of the datasets used in the paper is introduced in Section \ref{sec:dataset generation}. Sections \ref{sec:exploration} and \ref{sec:learning} contain an exploratory analysis of the dataset, normalization procedures and machine learning approaches to distinguish between Pareto and Dominated states. Finally, Section \ref{sec:conclusions} portraits some conclusive remarks and outlines future work.

\section{Literature} \label{sec:literature}

The body of literature dedicated to path finding in networks is massive, therefore, beside the already mentioned paper on the history of the shortest path \citep{schrijver2012history}, we limit ourselves to address the reader to \cite{Madkour2017}. 
Instead, in this section we focus on contributions to solve the (E)SPPRC, referring first to the surveys of \cite{Irnich2005} and \cite{Pugliese2013}.

Early approaches for SPPRC include pre-processing procedures and Branch \& Bound \citep{Christofides1989} and Lagrangian relaxations and enumeration of near-shortest paths \citep{Carlyle2008}. 
The most successful algorithms are based on \emph{Dynamic Programming labeling algorithms} \citep{Mehlhorn2000, Dumitrescu2003,  Desrochers1988} further improved by bi-directional search \citep{RighiniSalani2006} with dynamic balancing \citep{Tilk2017}, bucket based approaches \citep{PecinPessoa2017, Sadykov2020} and iterative relaxation techniques \citep{RighiniSalani2008}, \citep{Baldacci2011}, \citep{Martinelli2014}, \citep{Bulhoes2018} and \citep{Costa2021}.  For some applications, finding super-optimal non-elementary paths is relevant and in such cases efficient algorithms can be found in \citep{Irnich2006} and \citep{Desaulniers2008}

During the last years, Machine Learning (ML) techniques for the resolution of combinatorial optimization problems has become a popular trend. Its use can be classified, according to \citep{BENGIO2021}, as end-to-end, algorithm selection and configuration, and embedded in optimization algorithms.

Reviewing literature contributions for end-to-end ML approaches for CO goes beyond the scope of this work, but the interested reader can explore the trending topic of Neural Combinatorial Optimization starting from the survey by \cite{wu2024neuralcombinatorialoptimizationalgorithms}. A recent survey on the application of ML to Vehicle routing problems can be found in \cite{Reza2024}. For what concerns algorithm selection and configuration, ML is proficiently used in different ways. In \cite{Kruber2017}, the authors explore automatic Dantzig-Wolfe decompositions, predicting if a given problem is suitable for a successful decomposition approach and, in case of a positive result, what algorithm should provide the best one. \cite{NEURIPS2023_bcdec1c2} accelerates the solution of MILP problems by configuring the separators used to generate cutting planes. Similarly, \cite{NEURIPS2021_cb7c403a} proposes data-driven techniques for configuring the scheduling of heuristics procedures in exact MILP solvers.
Finally, machine learning techniques can be \emph{integrated in the optimization process}. In this setting, we can find works focused on integrating learning and branching \citep{Lodi2017} and studies where decision trees have been used to imitate strong branching policies in \cite{Alvarez2014}, resulting in fast approximations. In \cite{Morabit2021} instead machine learning is used to speed-up the re-optimization of the master problem during column generation, through a data driven selection of the most promising columns. Finally, the authors of \cite{Roman2018} focus on the pricing problem, exploiting knowledge gained from previous iterations to estimate upper bounds to reduce the solutions space in future iterations, through the means of generic online machine learning methods, providing considerable CPU time gains. 

To the best of our knowledge, very few works have tackled directly the integration of Machine Learning for Shortest Path Problems. Recently, \cite{morabit2023machine} studied how to speed-up shortest paths subproblems in vehicle routing and crew scheduling column generation algorithms. They exploit data driven classifiers trained on previous executions to perform arc selection, reducing effectively the size of the network and providing savings up to 40\% of computing time. In \cite{Xu2025}, the authors utilize reinforcement learning to select heuristics used to prune the network to speed up column generation. In \cite{morabit2024learning} instead, the authors propose a learned heuristic for the re-optimization of recurring routing problems, that present small differences in data, by predicting and fixing edges having high chances of being part of the solution at each iteration of the algorithm.
Indeed, we found no work that intersects machine learning with RCSPP to provide methodological improvements except that of \cite{ondei2024data} where authors attempt to forecast the feasibility of a partial path through data-driven techniques.

\section{Labeling algorithms for ESPPRC} \label{sec:algo}

The ESPPRC is defined on a graph $G(N,A)$, which we assume directed, composed of a node set $N$ and an arc set $A$. The problem asks to find a minimum cost elementary path, i.e. a finite sequence of consecutive arcs in which every node $n \in N$ appears at most once, from a source node $s \in N$ to a destination node $d \in N$. The cost is accumulated when traversing arcs along the path. We remark that no assumption are taken on the cost of the arcs and the graph may possess negative cost cycles. 

State of the art algorithms to solve the ESPPRC are based on a Dynamic Programming (DP) framework. 
Modern techniques include bi-directional search \citep{RighiniSalani2006}, \citep{Tilk2017}, bucket based approaches \citep{PecinPessoa2017, Sadykov2020} and advanced relaxation schemes \citep{RighiniSalani2008}, \citep{Baldacci2011} and their hybridizations \citep{Martinelli2014}.

In dynamic programming, a state associated with vertex $i \in N$ represents a partial path from the source node $s$ to the node $i$. Indeed, different states can be associated with the same node and they correspond to different partial paths. The DP algorithm iteratively extends states until no further extensions are possible and among all feasible states reaching the destination node $d$ the one with minimal cost represent the optimal solution to the ESPPRC.

Each state is encoded in a label, in bi-directional dynamic programming called {\it forward} and {\it backward} labels. Each of them stores the accumulated cost, the accumulated resources consumption and a binary vector that keeps track of the visited nodes in the partial path.
Formally a label $l^f_i$ representing a partial path ending at node $i$ is a tuple:
\begin{equation}
    l^f_i = (i, c, S, R)
\end{equation}
where $i$ is the last node visited in the partial path, $c$ is the accumulated cost, $S$ is a binary vector that keep tracks of the visited nodes in the partial path and $R$ is the so called resource vector that keeps track of the consumption of each resource. Bi-directional dynamic programming limits the extension of forward and backward states by selecting a monotone resource, called {\it critical resource}, and extending labels for which its consumption is less than a given threshold $T$. Such threshold can be adjusted dynamically in order to keep the dimensions of the forward and backward sets of labels as balanced as possible \citep{Tilk2017}.  

Regardless of the extension direction, given two labels $l_{1,i} = (i, c_1, S_1, R_1)$ and $l_{2,i} = (i, c_2, S_2, R_2)$ ending at the same node $i$, $l_1$ dominates $l_2$ iif:
\begin{align*}
c_1 \leq c_2\\
S_1 \subseteq S_2\\
R_{1,r} \leq R_{2,r} & \qquad \forall r \in R
\end{align*}
and at least one inequality is strict. Dominated labels can be safely discarded as they will not lead to an optimal solution.

Labels are managed in label \textit{pools}, normally associated to the nodes of the original graph or to the nodes of a specific data structure called bucket graph. Often, these label pools are implemented with priority queues so that a partial ordering between labels can be enforced.
Every label that is considered during the execution of the algorithm can be therefore classified as Generated and Dominated $(GD)$, that is, generated by the algorithm but discarded because a dominating label is already present in the label pool, Generated and Inserted $(GI)$, that is generated by the algorithm inserted in the label pool but later on dominated by a better label, Pareto $(P)$ that is a label generated by the algorithm, inserted in the label pool and not dominated in any pairwise dominance check. 

An effective strategy for addressing the ESPPRC involves iterative relaxation.  This approach aims at reducing the number of labels to be generated by dynamic programming by solving relaxations of the original problem and exploiting its optimal solution to drive the algorithm. The basic idea behind that scheme is that only a subset of the nodes of the network are relevant to compute the optimal solution without cycles. 
Therefore, the binary vector $S$ should be restricted only to those nodes. As a result, the domination criteria between labels is also relaxed and the algorithm generates less labels. At the same time, cycles on the nodes not considered in the set $S$ are not prohibited any longer.

The Decremental state space relaxation scheme (DSSR \citep{RighiniSalani2008}) relaxes a large portion of the state space, and iteratively adds nodes until a solution without cycles is found. Attempts to initialize the set $S$ have been explored in \citep{RighiniSalani2009}. Another relaxations scheme, called {\it ng-path} relaxation \citep{Baldacci2011}, consists of defining, for each node of the graph, a subset of nodes on which path elementarity is enforced. Hybridizations of the ng-path relaxation with the DSSR scheme are also possible \citep{Martinelli2014}, and, in particular, we focus on the Restricted DSSR (DSSR-R). DSSR-R stores a dedicated set for elementarity at each node and uses a DSSR based approach to iteratively forbid partial paths containing cycles. Namely, when a cycle is found, nodes that appear more than once are inserted in the sets of nodes appearing in the loop only.

%\begin{itemize}
%\item {\tt DSSR:} this is the standard DSSR scheme, where the set $S$ is initially empty and all nodes that are visited more than once in the optimal solution of each relaxation are added to the set $S$.
%\item {\tt DSSR-R:} stores a dedicated set for elementarity at each node and uses a DSSR based approach to iteratively forbid partial paths contaning cycles. Namely, when a cycle is found, nodes that appear more than once are inserted in the sets of nodes appearing in the loop only.
%\end{itemize}

In the following, we study how data gathered during the DSSR-R procedure can be exploited to determine whether a new state belongs to either GD, GI or P set.

In particular, we are interested in answering two research questions:
\begin{enumerate}
    \item is it possible to discriminate between Pareto and Dominated labels, with data driven techniques, by only looking at static properties of the labels?
    \item is it possible to determine whether a label about to be inserted in a label pool will eventually be dominated by a label that has not yet been generated, thus potentially allowing us to identify such cases in advance?
\end{enumerate}

To answer these questions, in the following, we specifically generate and analyze two datasets, G and I, and investigate the contexts in which different types of labels can be discriminated, considering both offline and online learning scenarios. %In offline learning we collect data related to states generated in the last iteration of the DSSR-R procedure and train a model that is used to classify states generated while solving different instances, while in online learning we collect data at an iteration $i$ of the DSSR-R procedure to train a model used to classify states generated at iterations $i+k$ of the DSSR-R procedure.

\section{Datasets Generation} \label{sec:dataset generation}

The first challenge of this work was the generation of a representative and diversified dataset, consisting of both Pareto and dominated labels for different instances, along with the careful design of meaningful features that could be computed in trivial time. Efficient feature computation is indeed crucial for this work, as recording states generated by dynamic programming algorithms can introduce overhead.
%, therefore making optimization essential.

\paragraph{Instances} More in detail, we generated our data starting from 41 problems of \cite{spprclibrepo}, an online repository of resource constrained shortest path problems. In this collection, all instances present complete graphs with cycles, having positive travel distances and prizes that can be collected (thus allowing the existence of cycles), and a single capacity constraint.
They are organized in different classes (A, B, E, G, M, P) that describe different instance generation approaches about capacity, node positioning and node clustering. Further details about how these sets were generated can be found in \cite{UCHOA2017}.

\paragraph{Labels generation}
We solved each problem through a modified version of PathWyse 1.0 \citep{Salani2024PathWyse}, an open source library for RCSPP, which employs a dynamic programming labeling algorithm and relaxation techniques, with the following settings: parallel bidirectional search, with an even fixed split of 50\% of the critical resource budget for each direction, DSSR-R (with two-cycle elimination), no NG-routes, node based extension strategy, Pareto join technique (see \cite{Salani2024ParetoJoin}). All runs were performed on a machine with an eight-core Intel i9-11900 @2.50 GHz, 32GB RAM and Kubuntu 24.04 OS. We used the same machine in every subsequent analysis. A one-hour time limit was set for the simulations, resulting in three unsolved instances. An additional instance was discarded during pre-processing operations, thus reducing the overall dataset to 41 instances. The final list of problems is reported in Section \ref{appendix:instances}, in the Appendix.

Since we used DSSR-R, for each algorithm iteration, we solved a different relaxation of the problem, and saved features for every label generated during optimization for both backward and forward search: to this end, we designed and implemented new data collection methods for recording and storing labels efficiently.

%We also conducted preliminary experiments on a dataset generated using the standard DSSR. In the following sections, most of the analyses and results apply to both versions of DSSR; we explicitly indicate in the text whenever this is not the case.

\paragraph{Features}

%Indeed, another challenge was to define, for each label, static properties, a set of characteristics that could describe a (partial) path while requiring no additional computational effort. 
We defined static properties, a set of characteristics that could describe a (partial) path while requiring no additional computational effort. In our setting, static features are accessed in constant time when a label is generated. This approach was a key element during feature design, offering the dual benefit of controlling overhead in dataset generation while potentially enabling low-latency machine learning applications for future algorithms. %The trade-off for this design choice is increased memory usage, which, given its abundance, was considered a reasonable compromise.

More specifically, for each label, we stored an identifier (iteration) that tracks the DSSR-R relaxation responsible for generating it, along with 23 additional features describing the path. These included the overall cost (objective), consumption of the critical resource (consumption\_critical), and tour-related features such as the length (tour\_length), the number of visited nodes, the number of unreachable nodes, and more. Additionally, we stored the full set of visited nodes used by the algorithm in the form of bitsets.
We also introduced an additional feature, named efficiency, which serves as a proxy for the quality of a label by accounting for both its cost and critical consumption. The definition and computation of this feature are detailed during the normalization process, as described in Section \ref{normalization}.

Features that detail the state of the network at the time a label was generated as well. These include the number of labels currently in the network for a forward or backward search (nlabels\_network), the number of dominated labels (nlabels\_dominated\_network), the number of extended labels (nlabels\_closed\_network), and the number of labels pending processing (nlabels\_open\_network). We also tracked similar information for the label pools, which were associated with each node by PathWyse.

Finally, we labeled each entry (i.e. partial path) under the {\tt dominated} feature, marking labels from sets GD and GI as Dominated (1), and labels from set P as Pareto (0). These values are then used as classes in the next analyses, allowing us to distinguish between dominated and non-dominated labels.

The full list of features can be found in the Appendix, in Section \ref{appendix:features}. After preliminary experiments, we performed feature selection by discarding a small subset of features that could either not be used directly, such as the full bitsets, or identifiers (for the run or the nodes). We also performed Principal Component Analysis in Section \ref{sec:exploration}. Overall, we marked with an asterisk the features in the Appendix that were selected for analysis with machine learning models in Section \ref{sec:learning}.

\paragraph{Datasets} We consolidated forward and backward labels of all instances, from every DSSR-R iteration, into two primary datasets:
\begin{itemize}
    \item \textbf{Dataset G} contains every label that has been generated by the labeling algorithm (P, GD, GI).
    \item \textbf{Dataset I} includes only labels that have been inserted in data structures (P, GI). Therefore, Dataset I is a subset of Dataset G.

\end{itemize}

A summary of the datasets is reported in Table \ref{tab:summary}: for both datasets and every instance class, we report the number of problems, the overall number of generated labels, the percentage of Pareto ones, the average number of DSSR-R iterations required to solve an instance and the average number of labels available at each DSSR-R iteration.

\paragraph{Dataset G}

This dataset includes all the labels that are generated by the algorithm: it contains all Pareto labels (P) and all dominated labels that are either discarded before insertion (Generated and Dominated, GD) or afterwards (Generated and Inserted, GI), when a new, more promising label is inserted in the same label pool. This means, that all dynamic programming states, even suboptimal ones, are recorded for both forward and backward search.

\begin{figure} [htb]
\begin{subfigure}{.49\textwidth}
  \centering
  \includegraphics[width=\linewidth]{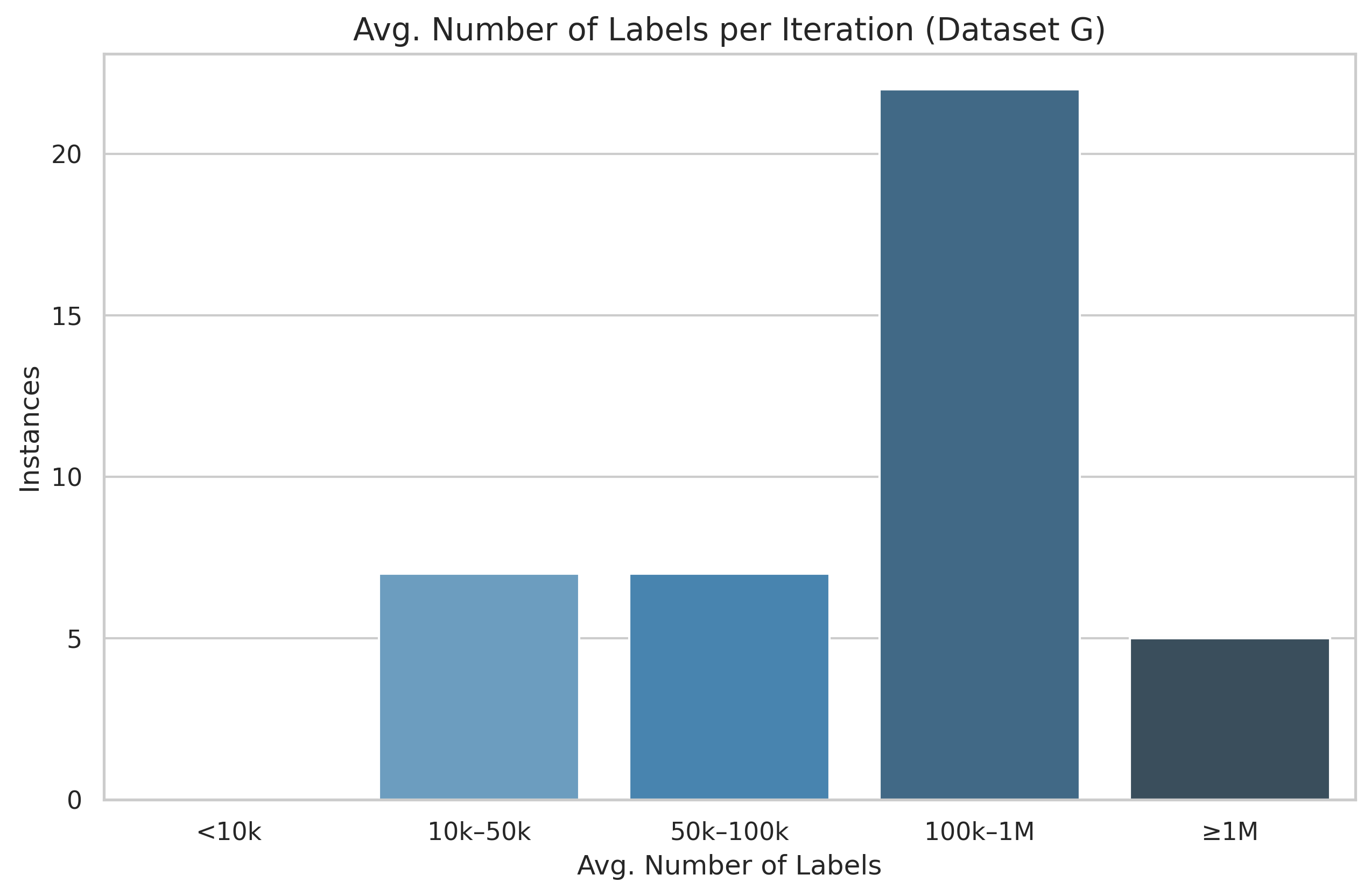}
  \caption{Dataset G}
  \label{fig:avglabG}
\end{subfigure}%
\begin{subfigure}{.49\textwidth}
  \centering
  \includegraphics[width=\linewidth]{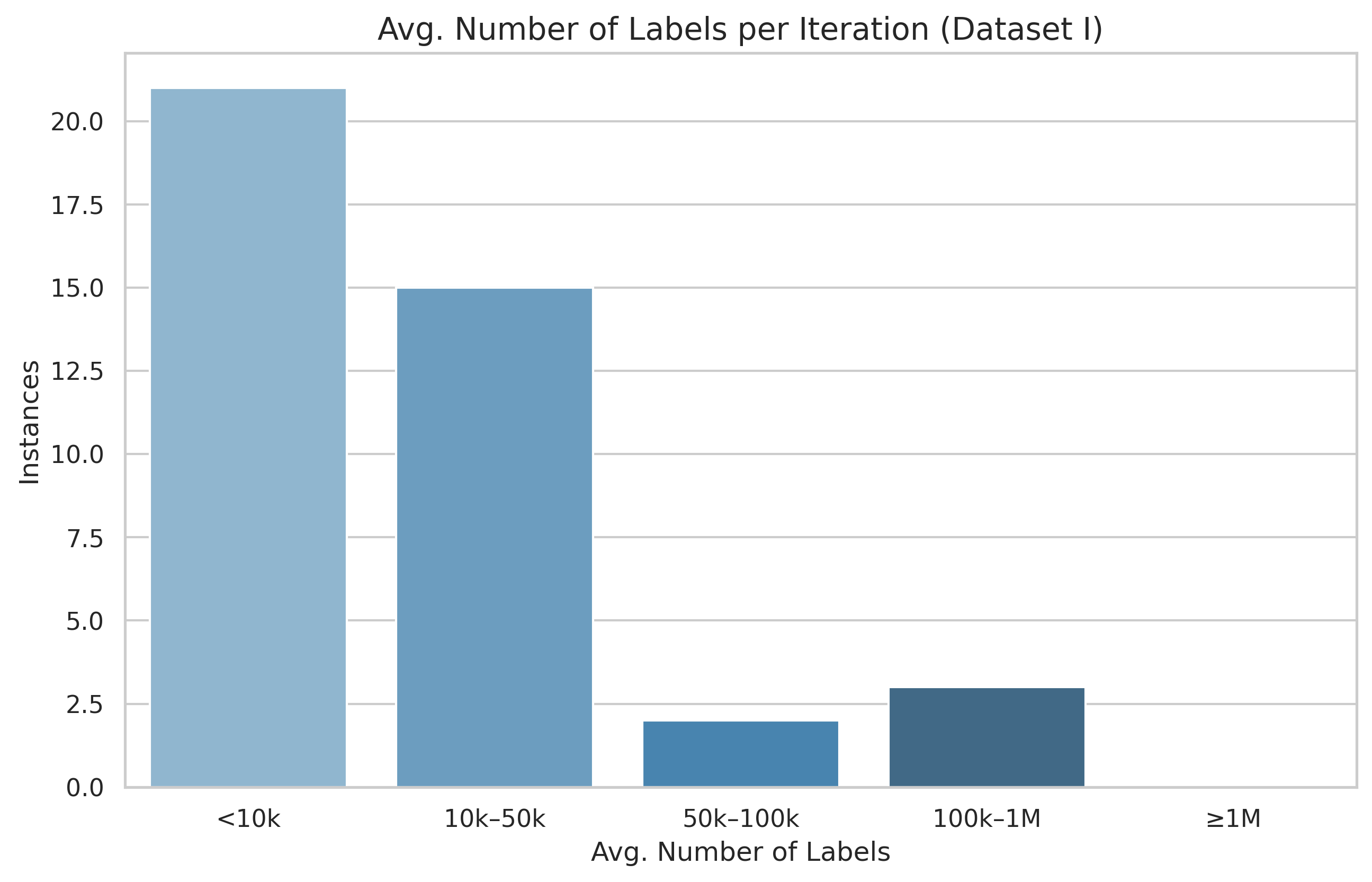}
  \caption{Dataset I}
  \label{fig:avglabI}
\end{subfigure}
\caption{Average number of labels per iteration, for each instance.}
\label{fig: profiling}
\end{figure}

Therefore, as reported in Table \ref{tab:summary}, this dataset consists of 41 instances, for about 1 billion records. Problems were solved, on average, in approximately 15 iterations of DSSR-R, with around 930,000 labels available per iteration. Indeed, the majority of instances exhibited more than 100,000 labels per iteration: this is summarized by the histogram shown in Figure \ref{fig:avglabG}. We note that instances of class M have one order of magnitude more records compared to the other classes.
We also observe that the dataset is clearly unbalanced, as the negative class Pareto represents only about 0.89\% of total records of the dataset. Indeed, out of all the possible states, only a very small subset of labels is dominating, confirming the importance of strong extension and dominance rules.

We remark that, overall, generating Dataset G required, on average, 54\% more time with respect to a standard PathWyse optimization run. This was expected, given the high number of records, and overhead can be reduced if data storage on disk is not necessary. Additional memory and efficiency strategies can be adopted if one is interested in storing only the labels from a single iteration.
%We designed this dataset to focus on a broad, fundamental research question, that is, understanding if it is possible to discriminate between Pareto and dominated labels, with data driven techniques, by only looking at static properties of the labels.

\paragraph{Dataset I}
This dataset collects instead, only labels that passed successfully the insertion step, that is, at generation time, there were no other labels dominating them. Therefore, it contains Pareto labels (P) and labels that were inserted in data structures but discarded after a more promising label was generated at a later point of the algorithm (Generated and Inserted, GI). 

The dataset contains data for all 41 instances, for a total of 32 million records. In general, there are two order of magnitude fewer labels compared to Dataset G. Only about 30,000 records are available for each DSSR-R iteration. This trend is illustrated in Figure \ref{fig:avglabI} where only few instances have more than 50,000 records per iteration. This pattern is especially prominent for instances of class E, whereas problems of class M have the largest sets of labels.

Unlike Dataset G, this dataset is more balanced: on average, Pareto labels account for over 30\% of the total. Given the smaller number of labels, the generation of Dataset I is on average only 25\% slower compared to the standard PathWyse solution.

%We designed this dataset to address a more specific research question: determining, based solely on static properties, whether a label about to be inserted into a data structure will eventually be dominated by a more promising label that has not yet been generated, thus potentially allowing us to identify such cases in advance.

\begin{table}[htb]
\small

\begin{tabularx}{\textwidth}{X*{7}{r}}\toprule

Dataset	&Class	&Inst.	&Labels	&Pareto	&DSSR-R It.	&Labels per It.\\
\midrule
\multirow{4}{*}{G}
&A	&10	&28519876	&2.80\%	&13.50	&210351.84\\
&B	&7	&26868310	&3.11\%	&15.86	&227244.19\\
&E	&6	&8007201	&2.91\%	&10.33	&205046.04\\
&M	&3	&1051067897	&0.69\%	&40.00	&8461281.09\\
&P	&15	&115294847	&1.57\%	&13.07	&522248.51\\

        \addlinespace
	Overall &	&41	&1229758131	&0.89\%	&15.22	&930294.51\\
        \hline
        \addlinespace
\multirow{5}{*}{I}
&A	&10	&1824155	&43.70\%	&13.50	&13480.59\\
&B	&7	&2060712	&40.55\%	&15.86	&17848.01\\
&E	&6	&475798	&48.90\%	&10.33	&10331.22\\
&M	&3	&23303860	&31.07\%	&40.00	&188299.21\\
&P	&15	&4704911	&38.50\%	&13.07	&21588.57\\

 	 \addlinespace
	Overall &	&41	&32369436	&33.73\%	&15.22	&29523.30\\
\bottomrule
\end{tabularx}
\caption[Summary] {Summary of Dataset G and Dataset I.  For each problem class, number of instances, number of labels, percentage of Pareto labels, average number of DSSR-R iterations and average number of labels per iteration.} \label{tab:summary}
\end{table}

\section{Exploratory analysis of the dataset} \label{sec:exploration}

In this section, we present an exploratory investigation of the datasets, covering key feature, instance and correlation analysis. These experiments are primarily conducted on a per-instance basis, to optimize memory usage, and on the last iteration of the algorithm, where solving the relaxation yields the optimal solution, although we also include some study across iterations (DSSR-R iterations). Finally, we design and detail normalization procedures at the end of this section.
All the analyses were performed in Python 3.10.
%why last iteration?

\paragraph{Analysis 1: Pareto vs Dominated labels, by instance}

\begin{center}
\begin{figure}
\begin{subfigure}{\textwidth}
  \centering
  \includegraphics[width=0.6\linewidth]{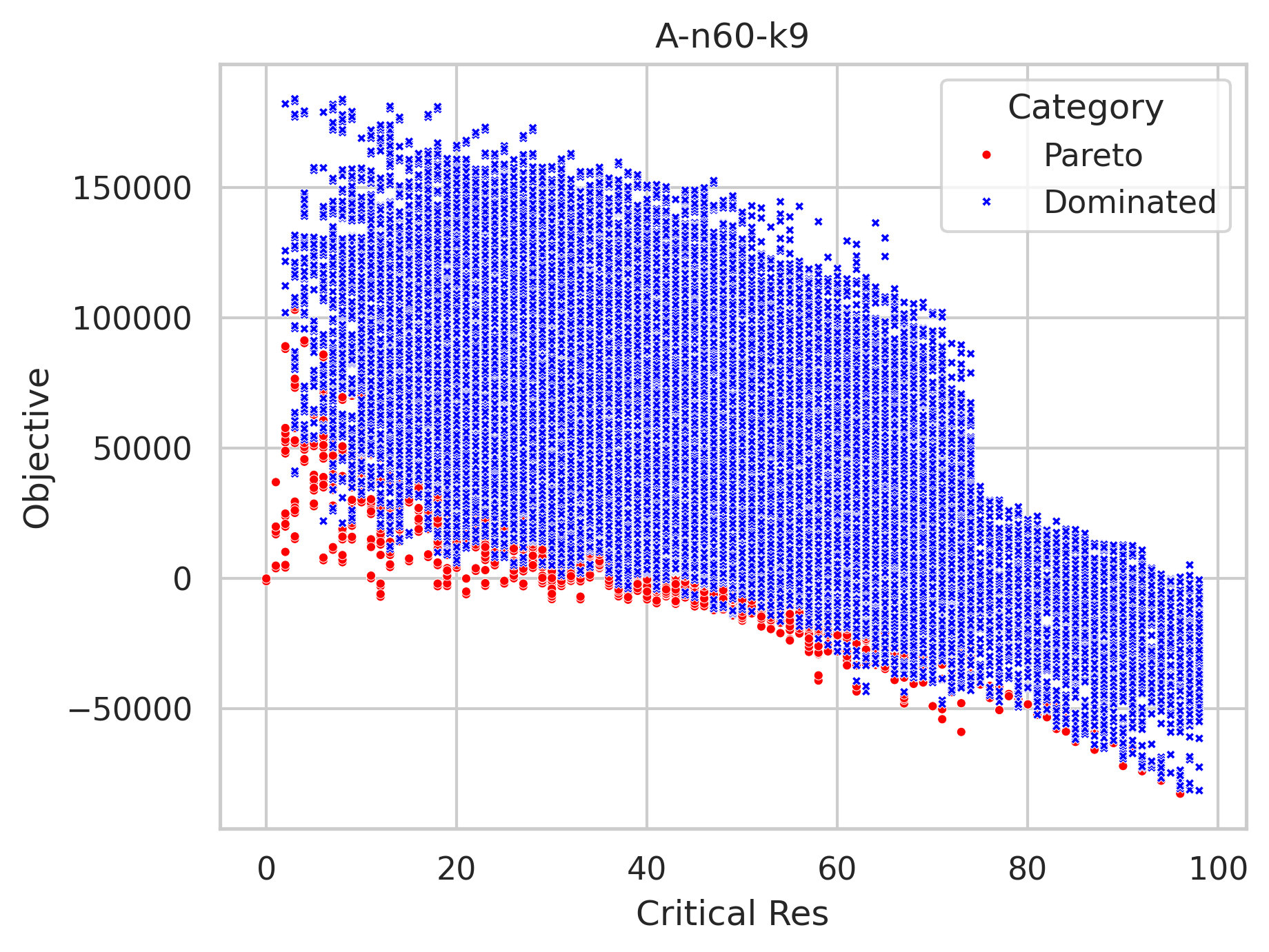}
  \caption{}
  \label{fig:scatterG1}
\end{subfigure}
\begin{subfigure}{\textwidth}
  \centering
  \includegraphics[width=0.6\linewidth]{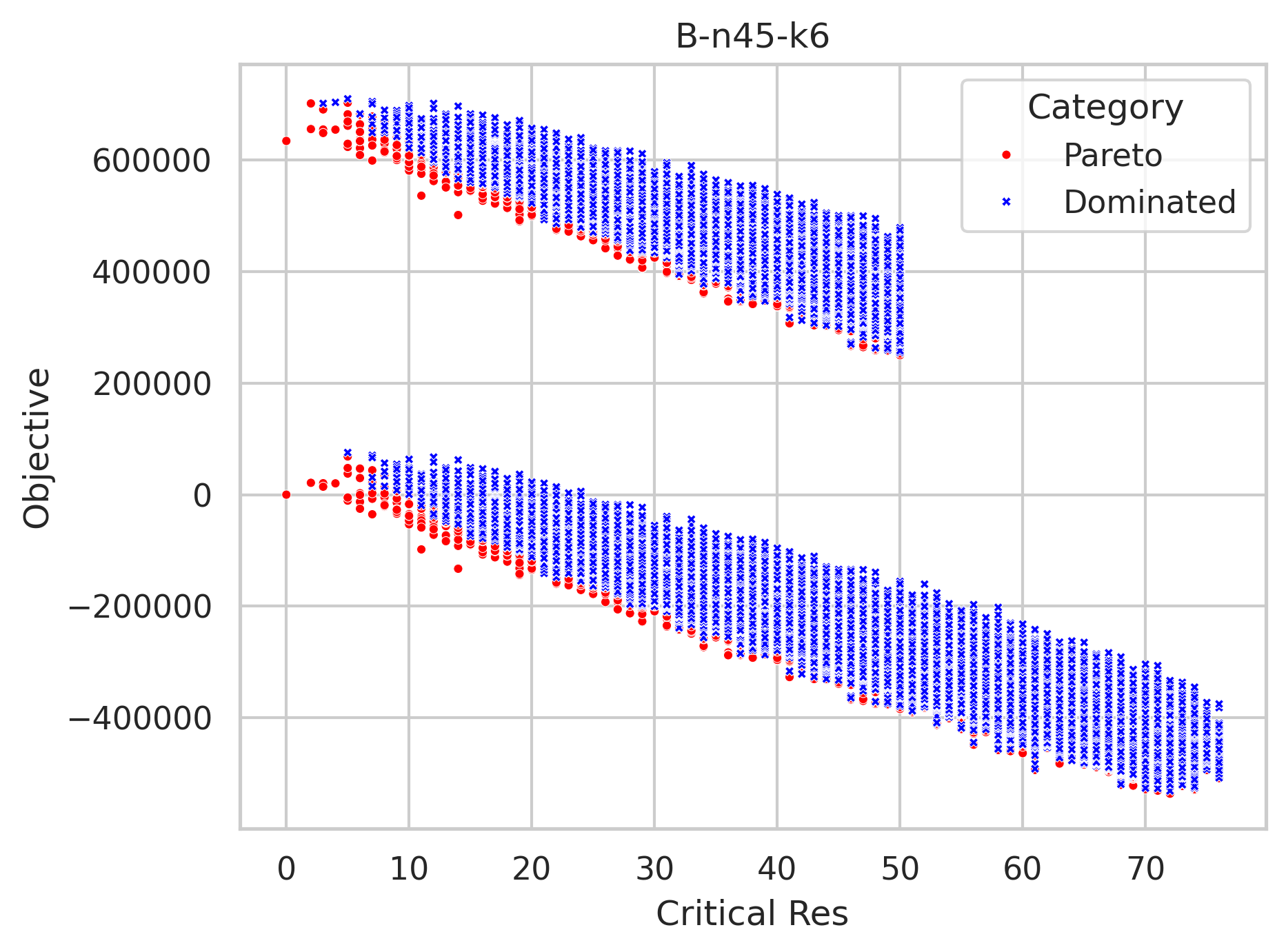}
  \caption{}
  \label{fig:scatterG2}
\end{subfigure}
\begin{subfigure}{\textwidth}
  \centering
  \includegraphics[width=0.6\linewidth]{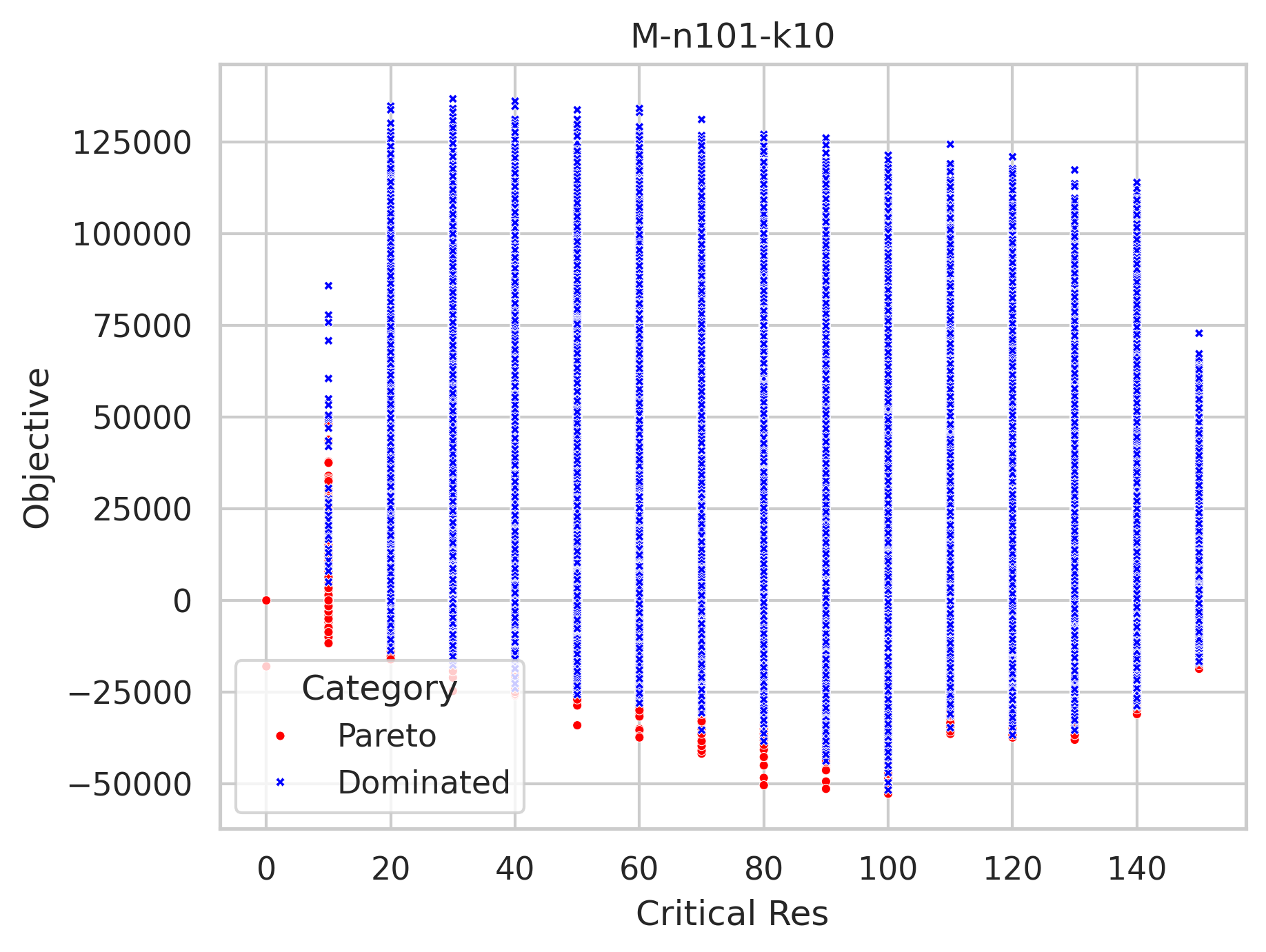}
  \caption{}
  \label{fig:scatterG3}
\end{subfigure}
\caption{Dataset G. Scatter plot of Label objective (y-axis) against critical resource (x-axis). Pareto labels are marked in red, while dominated ones in blue.}
\label{fig:scatterG}
\end{figure}
\end{center}

\begin{center}
\begin{figure}
\begin{subfigure}{\textwidth}
  \centering
  \includegraphics[width=0.6\linewidth]{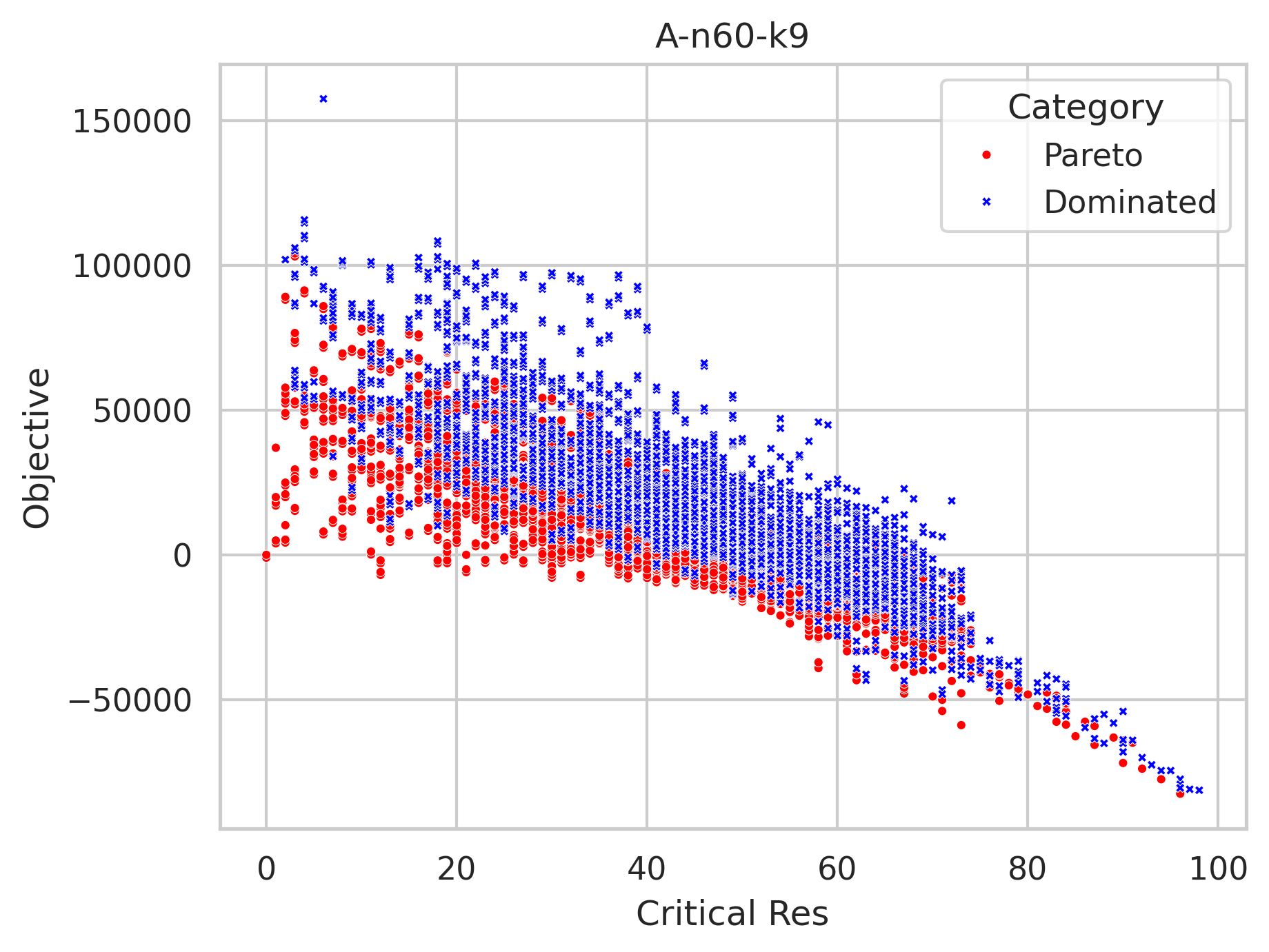}
  \caption{}
  \label{fig:scatterI1}
\end{subfigure}
\begin{subfigure}{\textwidth}
  \centering
  \includegraphics[width=0.6\linewidth]{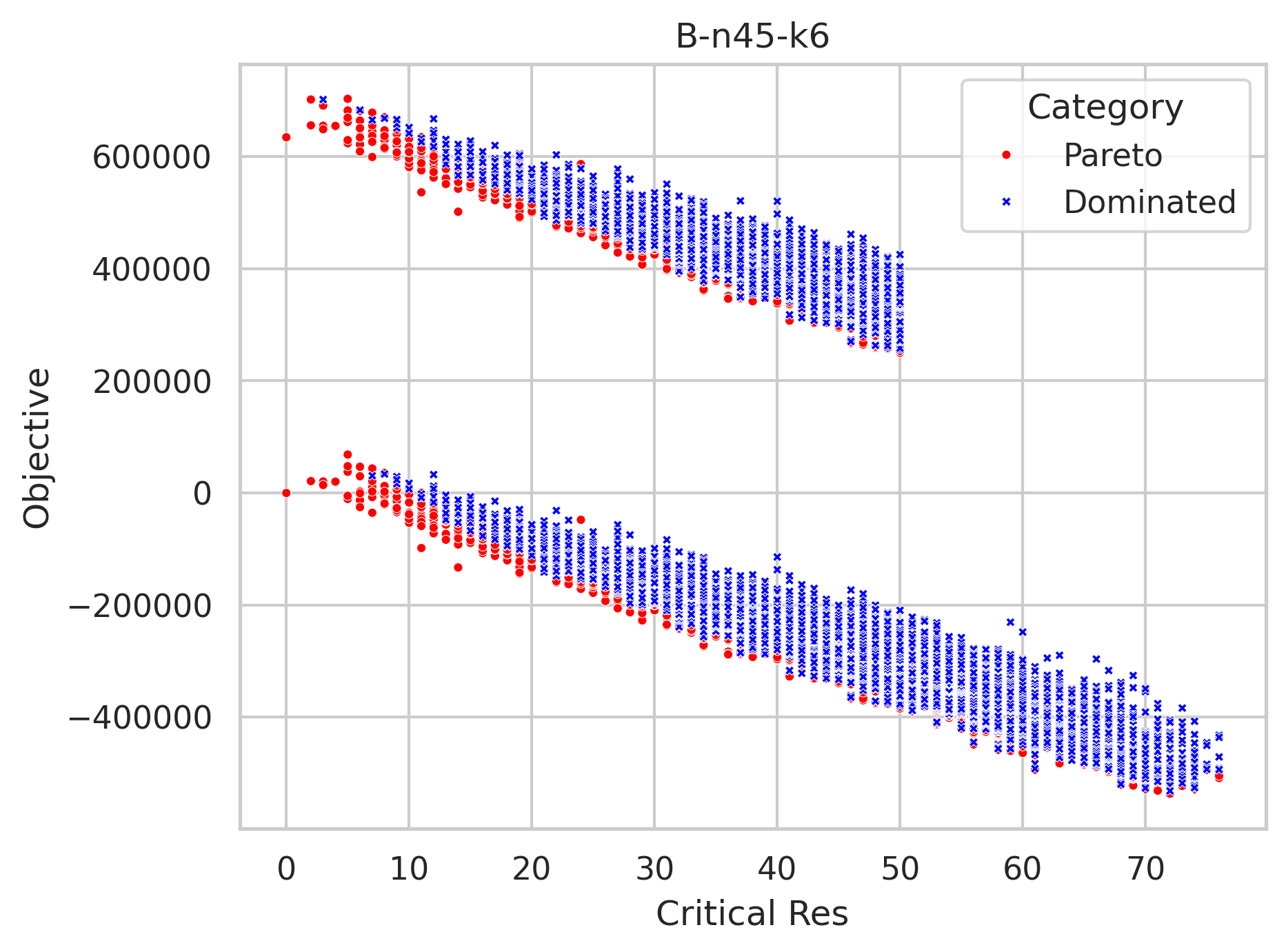}
  \caption{}
  \label{fig:scatterI2}
\end{subfigure}
\begin{subfigure}{\textwidth}
  \centering
  \includegraphics[width=0.6\linewidth]{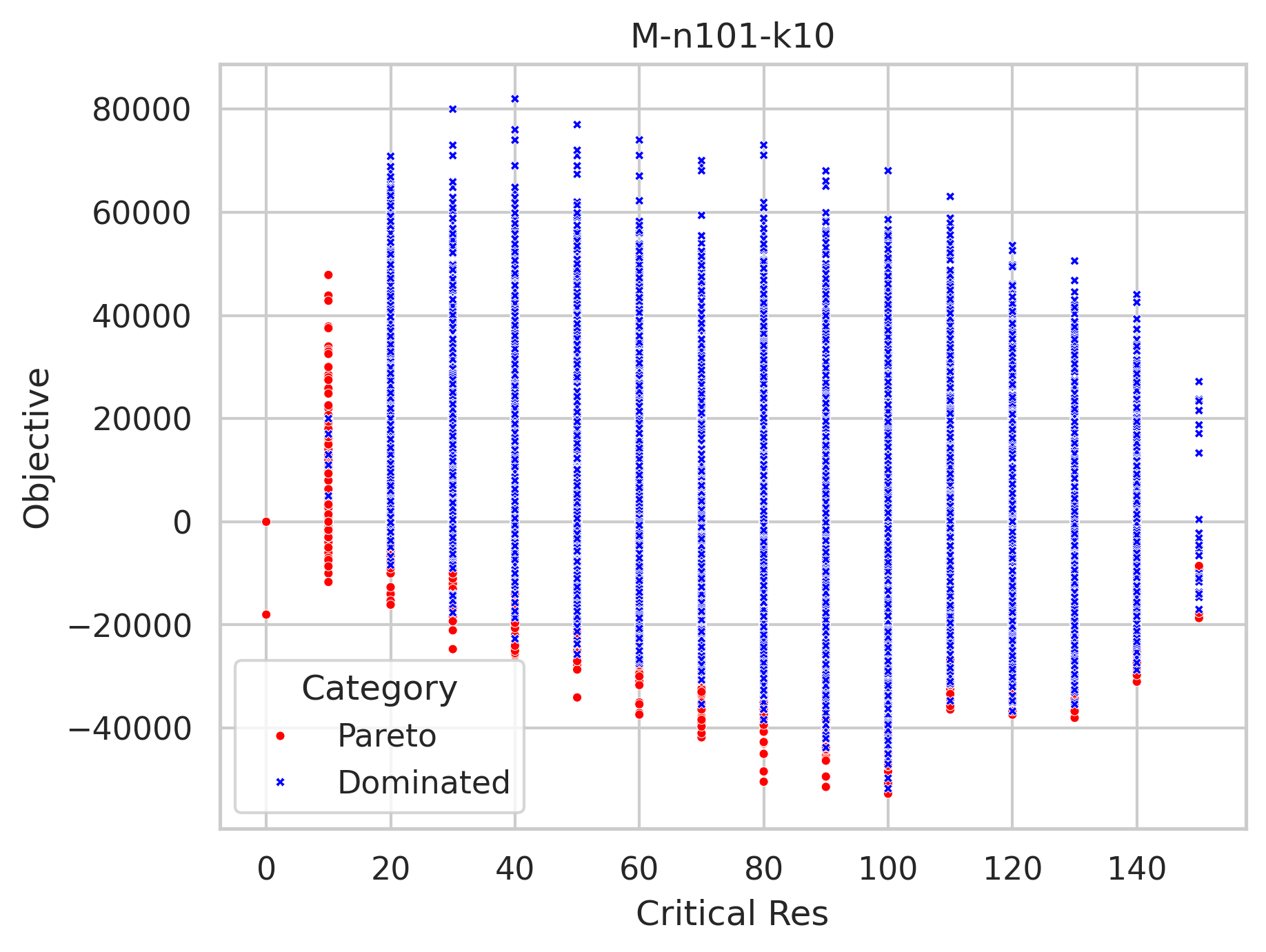}
  \caption{}
  \label{fig:scatterI3}
\end{subfigure}
\caption{Dataset I. Scatter plot of Label objective (y-axis) against critical resource (x-axis). Pareto labels are marked in red, while dominated ones in blue.}
\label{fig:scatterI}
\end{figure}
\end{center}

In the following, we summarize visually the distribution of labels for representative instances of Dataset G (resp. Dataset I) in Figure \ref{fig:scatterG} (resp. Figure \ref{fig:scatterI}). Each label is represented with its objective value on the y-axis and its critical resource consumption on the x-axis. Pareto labels are marked in red, while Dominated labels are shown in blue.
For Dataset G, we observe that the few existing Pareto labels generally have low objective values; as expected, when the cost increases, most labels are Dominated. An exception to this behavior is shown in Figure \ref{fig:scatterG2} (and two others), where the two searches (forward and backward) produce distinct sets of labels. This difference arises from the varying prize distributions around the source and destination nodes. For example, nodes near the source may have higher prizes, while those closer to the destination may exhibit strictly lower prizes.
In the following, we refer to these problems as {\tt split} instances.
Still, for this dataset, 37 problems out 41 were similar to \ref{fig:scatterG1}.
Finally, a single instance, {\tt M-n101-k10}, shown in Figure \ref{fig:scatterG3} presents instead a stack-like structure, since resource consumption is represented by discrete multipliers of 10. 

For dataset I, the remaining Dominated labels are much more intermingled with the Pareto labels, suggesting a greater overlap between the two sets. While most instances (35) are similar to Figure \ref{fig:scatterI1}, some (5) reveal the same split distribution (see Figure \ref{fig:scatterI2}) as reported for Dataset G. Finally, as before, only {\tt M-n101-k10} shows a stack-like structure. 

Overall, we observe that, as also highlighted in Table \ref{tab:summary}, most labels generated by dynamic programming algorithms are suboptimal (approximately, 99\% of the total). Given the small size of Dataset I, it is notable that state-of-the-art algorithms already employ robust rules to efficiently discard dominated labels before insertion, leaving mostly promising ones. While this is expected, we believe there is still potential for improvement by accelerating insertion checks and further reducing the number of insertions through the right selection of candidate labels for extension. Indeed, choosing the right candidate remains a critical factor in controlling the combinatorial explosion of labels.

\begin{figure}
\begin{subfigure}{.49\textwidth}
  \centering
  \includegraphics[width=\linewidth]{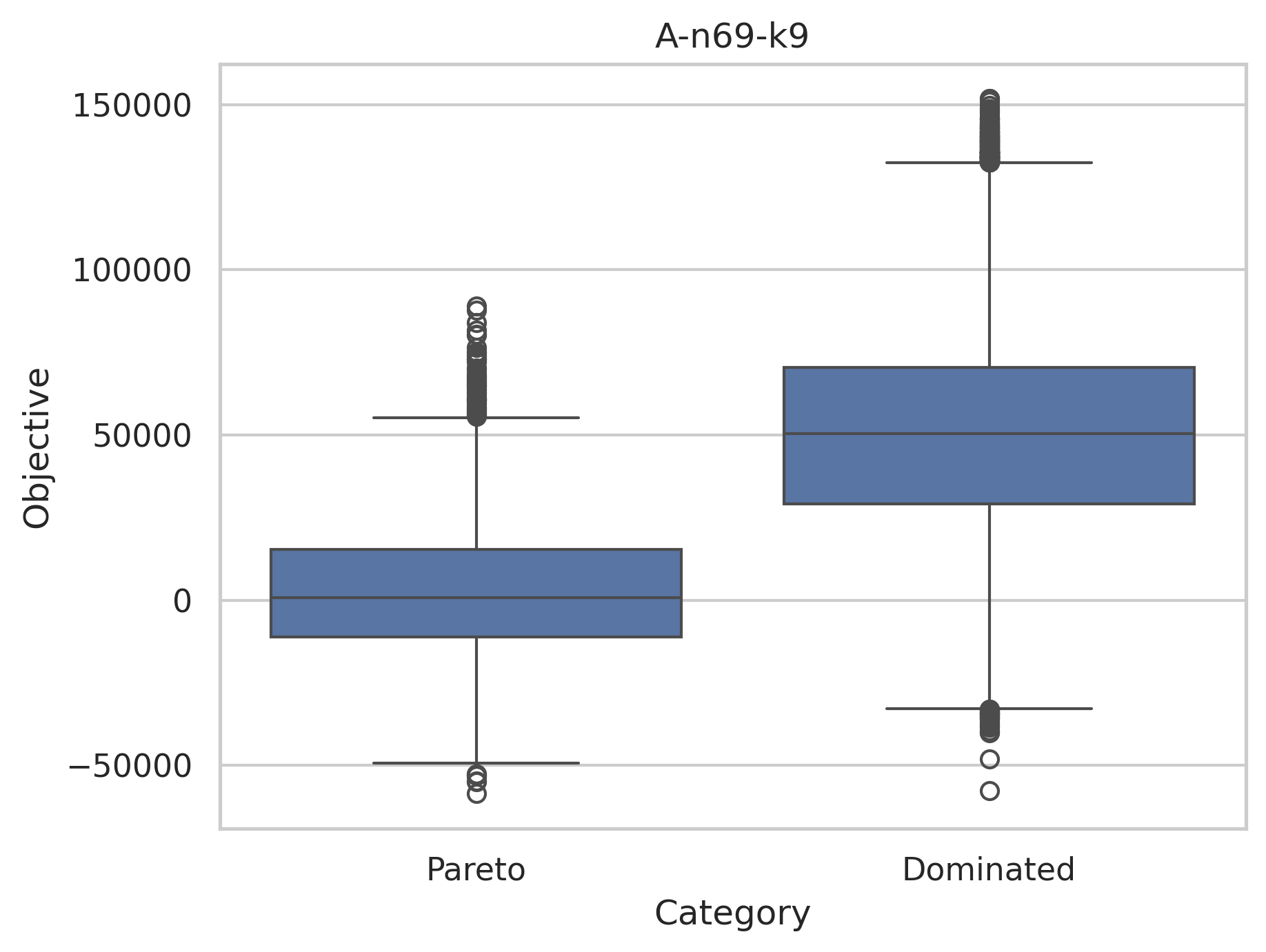}
  \caption{}
  \label{fig:boxobjGA}
\end{subfigure}%
\begin{subfigure}{.49\textwidth}
  \centering
  \includegraphics[width=\linewidth]{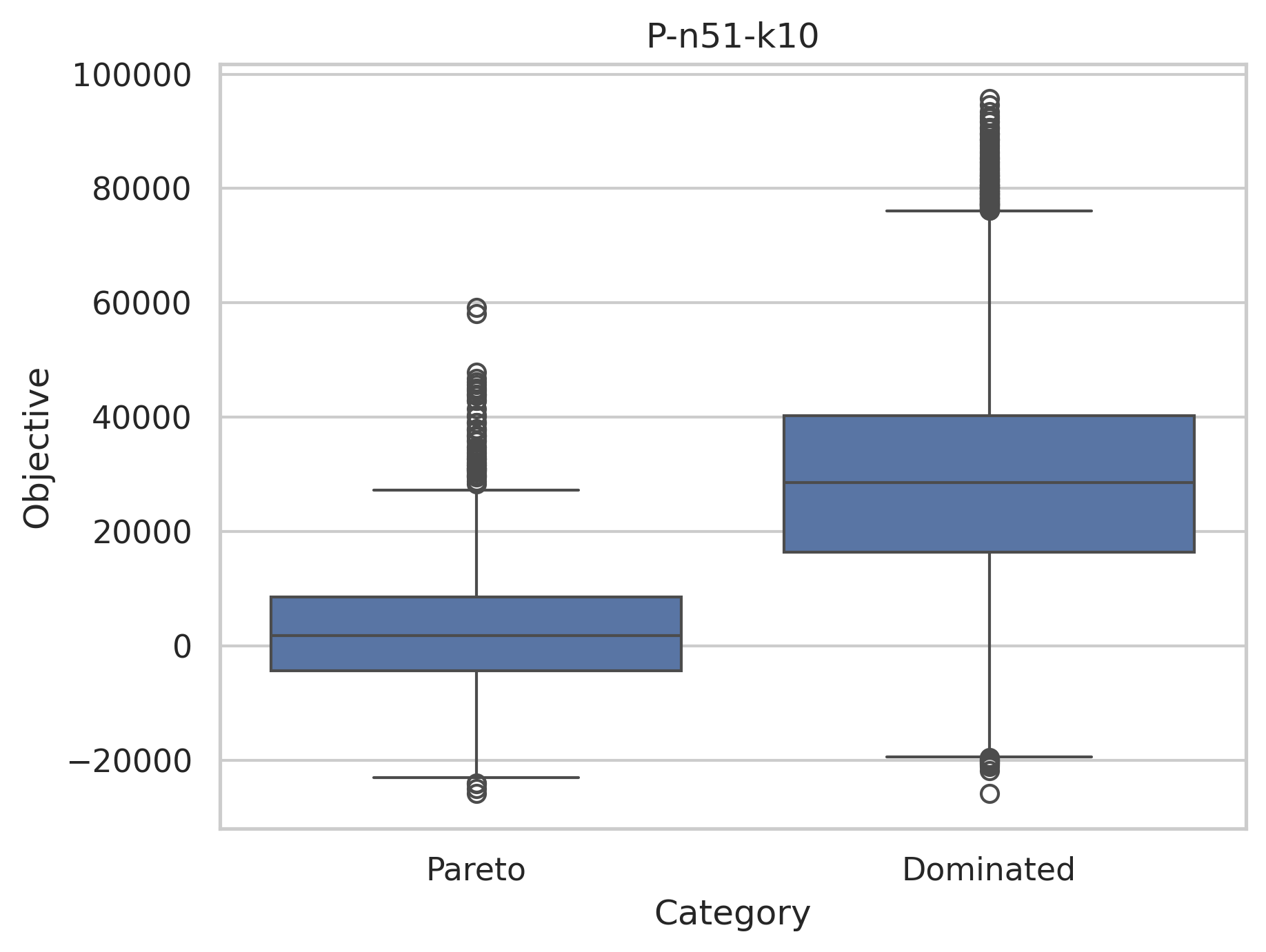}
  \caption{}
  \label{fig:boxobjGP}
\end{subfigure}
\caption{Dataset G. Boxplot of objective for Pareto and dominated labels.}
\label{fig: boxplotobjG}
\end{figure}

\begin{figure}
\begin{subfigure}{.49\textwidth}
  \centering
  \includegraphics[width=\linewidth]{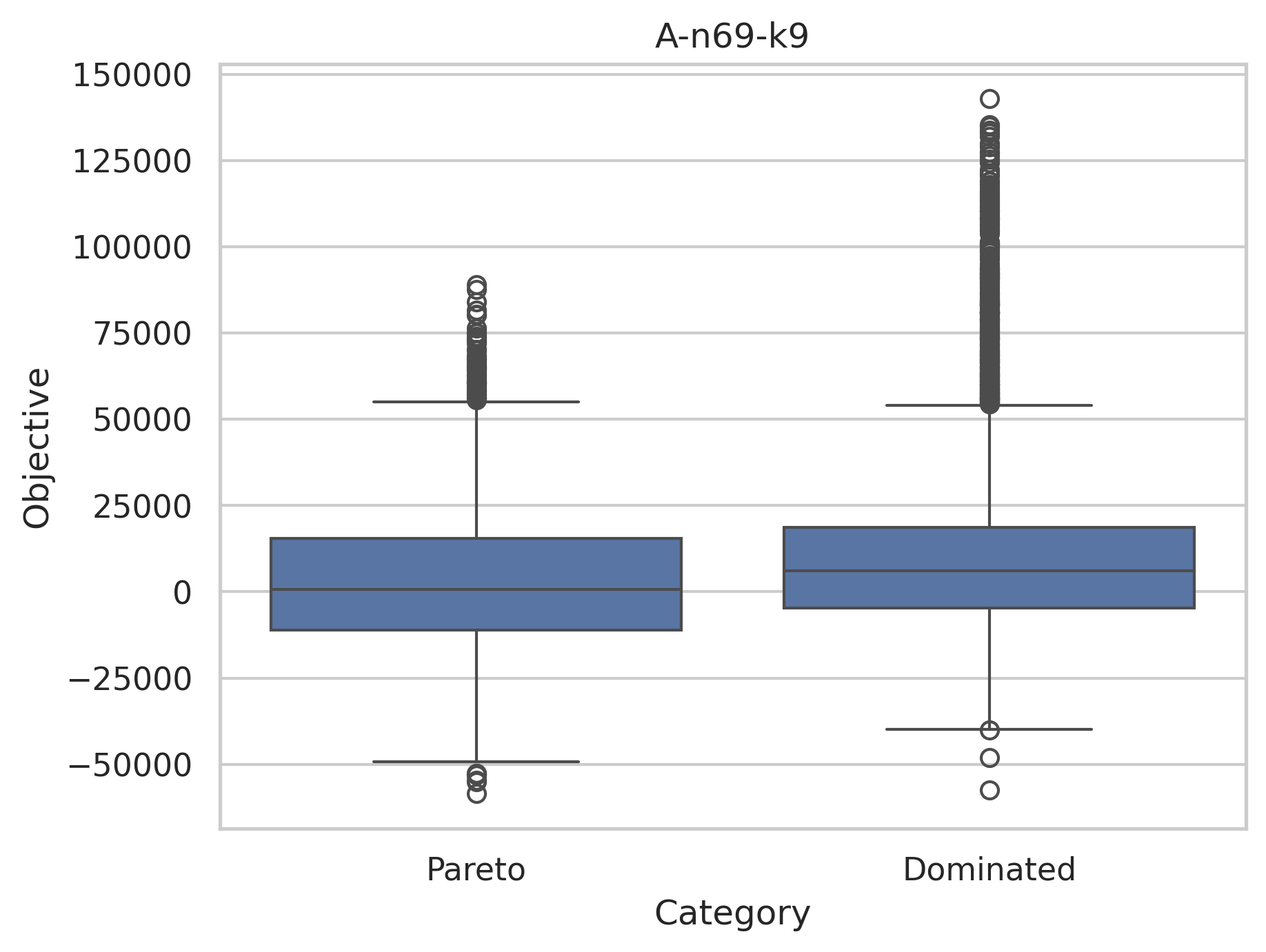}
  \caption{}
  \label{fig:boxobjIA}
\end{subfigure}%
\begin{subfigure}{.49\textwidth}
  \centering
  \includegraphics[width=\linewidth]{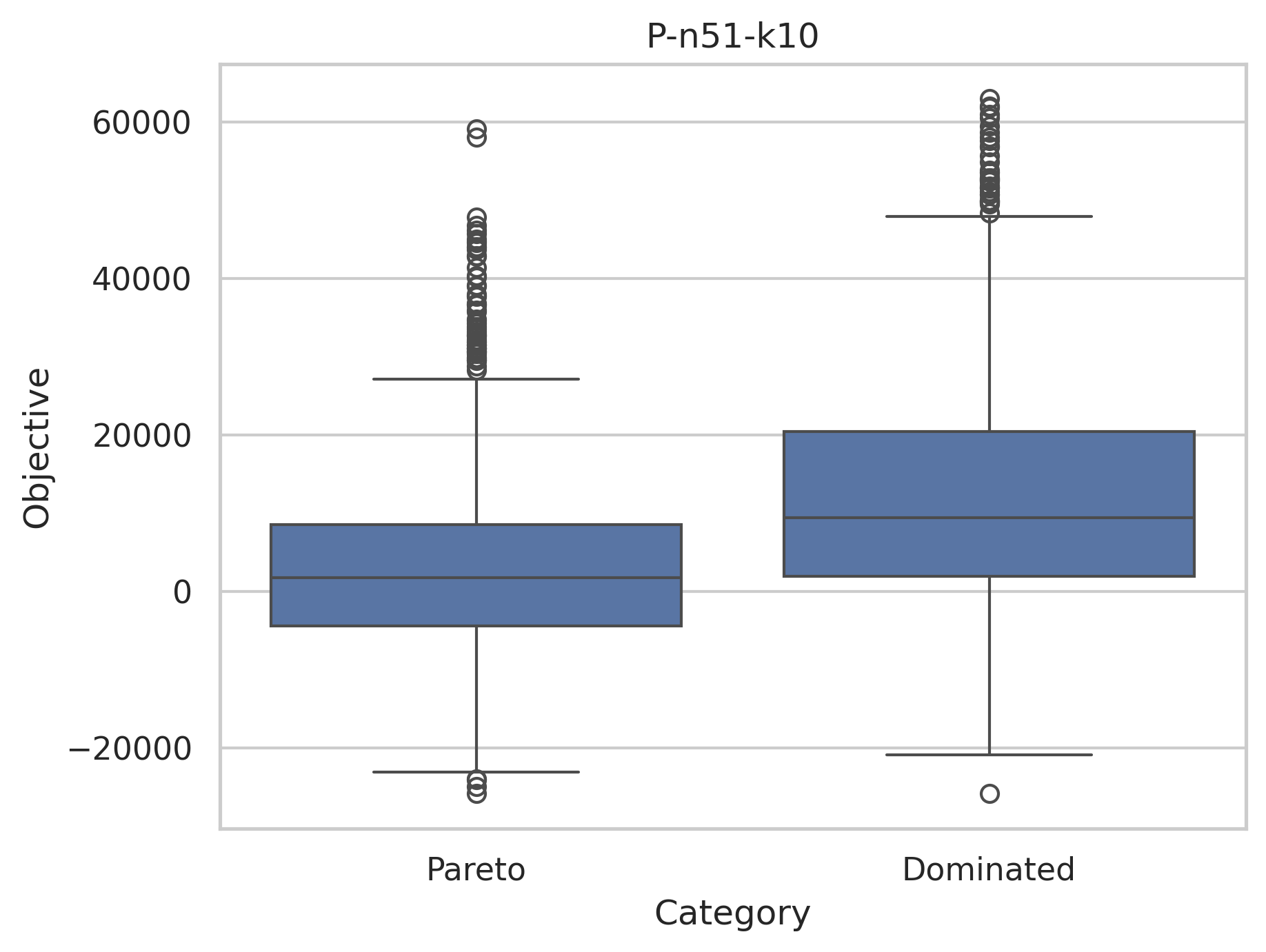}
  \caption{}
  \label{fig:boxobjIP}
\end{subfigure}
\caption{Dataset I. Boxplot of objective for Pareto and dominated labels.}
\label{fig: boxplotobjI}
\end{figure}

\paragraph{Analysis 2: Features distribution}
We conducted a more detailed analysis of the feature distributions utilized by the dominance functions during the algorithm's execution. In particular, we focused on the objective distributions for Pareto and Dominated labels in both Dataset G and Dataset I. The corresponding boxplots, for selected instances, are presented in Figure \ref{fig: boxplotobjG} and Figure \ref{fig: boxplotobjI}, respectively.
For Dataset G, objective distribution between the two type of labels is very different: in fact, we observe that the third quartile (Q3) of the Pareto-optimal ones is lower than the first quartile (Q1) of the dominated labels. This is true for SPPRCBLIB, for all instances but {\tt split} ones, due to forward and backward labels having different cost distribution.
Indeed, dominated labels have a very high objective, and simple rule-based thresholds might be able to discriminate most labels between the two sets.

This is not consistent with Dataset I, where Pareto-optimal and dominated labels show similar objective distribution. This suggests that insertion is efficient, leaving only labels that are quite similar to Pareto-optimal ones: in this case, the objective cost alone is not sufficient as a proxy to describe the quality of a label.

\begin{figure}
\begin{subfigure}{.49\textwidth}
  \centering
  \includegraphics[width=\linewidth]{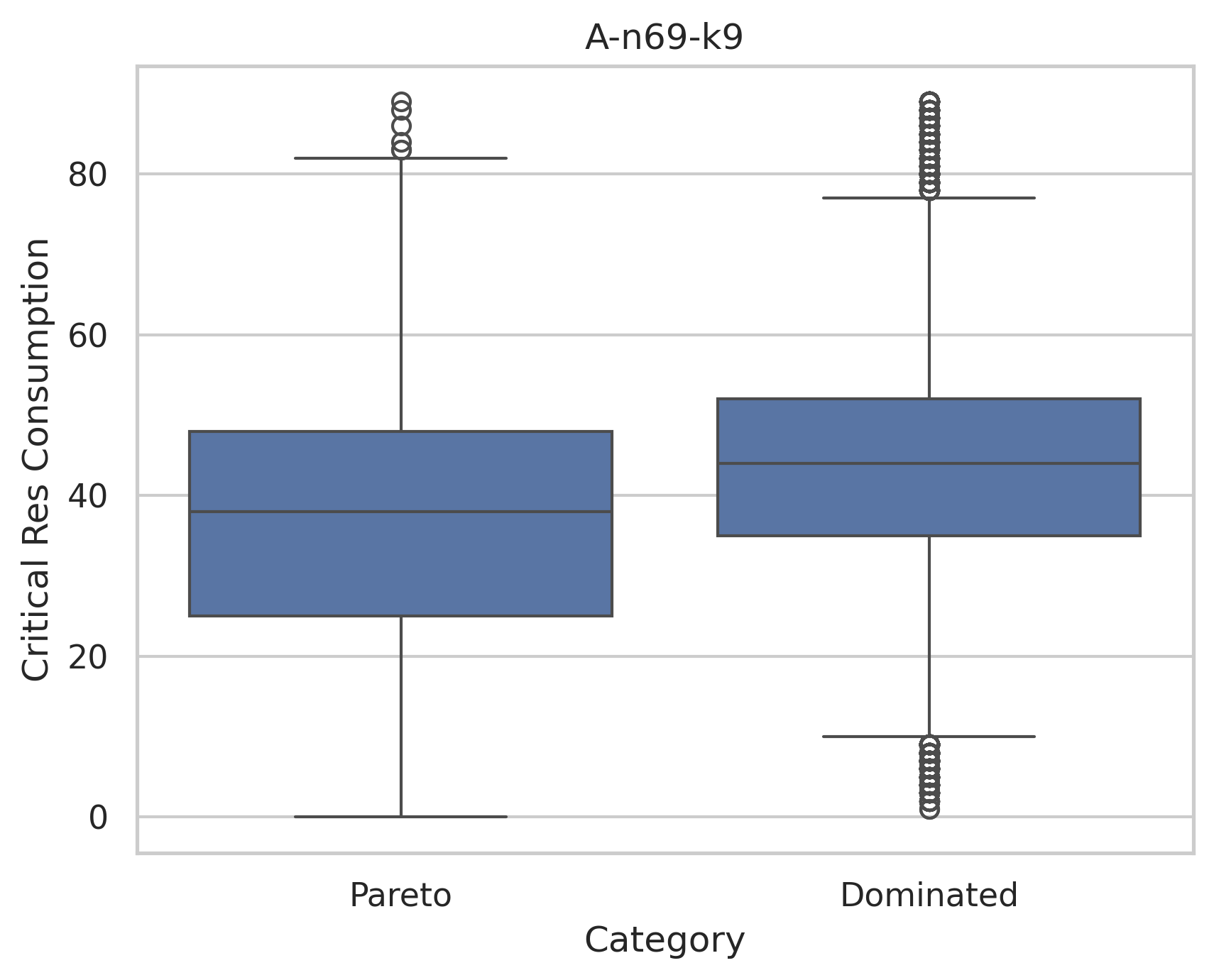}
  \caption{Dataset G}
  \label{fig:boxcritGA}
\end{subfigure}%
\begin{subfigure}{.49\textwidth}
  \centering
  \includegraphics[width=\linewidth]{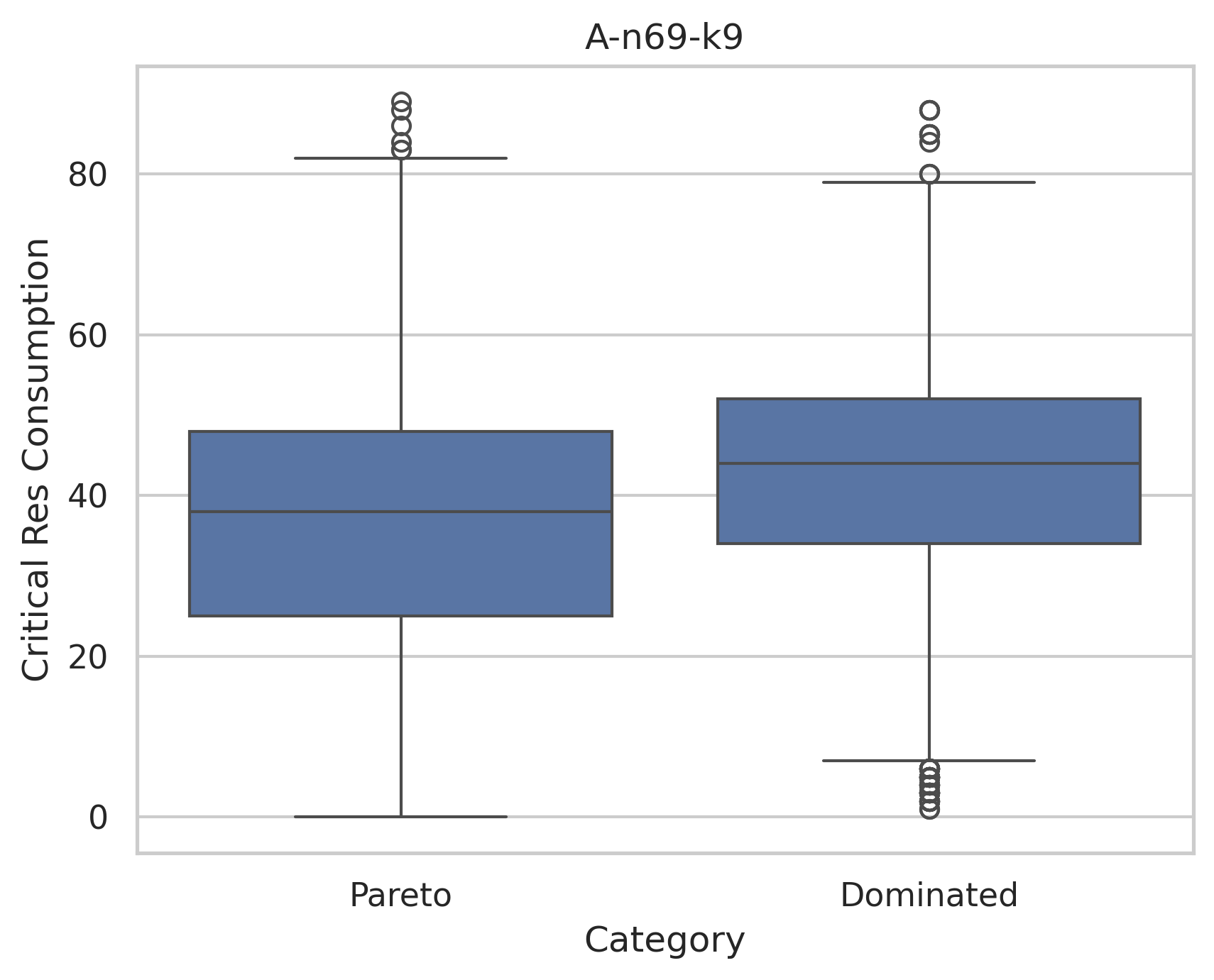}
  \caption{Dataset I}
  \label{fig:boxcritGP}
\end{subfigure}
\caption{Boxplot of critical resource consumption for Pareto and dominated labels.}
\label{fig: boxplotcrit}
\end{figure}

We also examined other features, including critical resource consumption, the number of visited nodes, and tour length. However, in every case, Pareto-optimal and dominated labels showed similar distributions, making it challenging to discriminate between the two sets. This observation applies to both datasets and is illustrated in Figure \ref{fig: boxplotcrit} for critical resource consumption.

\paragraph{Analysis 3: Correlation analysis and Principal Components Analysis}
For each instance, we also measured the linear correlation between each feature and the binary label type to assess how well each feature separates the Pareto (0) and Dominated (1) classes.

Results are reported in Table \ref{tab:corr} for both Dataset G and I. For each feature, we report the average, minimum, maximum and average correlation across all instances. We remark that values closer to $1$ or $-1$ are desirable, as they indicate stronger positive or negative correlations, respectively.

Indeed, we found that no single feature exhibits a strong correlation with the label type. Specifically, for Dataset G, most features show an average correlation close to 0, albeit the objective is the only feature with weak positive correlation. The number of unreachable and visited nodes show instead weak negative correlation.
Still, objective, tour length, number of unreachable nodes and number of visited nodes show a weak to moderate relationship with the label type on a limited number of instances.

Although the results for Dataset I are similar, we observe that, unlike Dataset G, the objective correlation is lower. More importantly, we note a moderate negative correlation, on average, with the features that describe the state of the network. Specifically, the values suggest that as the number of labels increases and the label pools become richer in Pareto labels, the algorithm improves at performing insertion, leading to fewer dominated labels being inserted into the pools.

We also tested Principal Component Analysis (PCA) and, based on the eigenvalues criterion, identified 10 significant components. The average absolute score of the components most strongly correlated with the label type was 0.2 for Dataset G and 0.33 for Dataset I. Given the dataset’s relatively modest dimensionality and the lack of improvement in correlation scores, we conducted the subsequent experiments using the full set of selected features.

Overall, this analysis suggests that the task to discriminate between Pareto and dominated labels is indeed a non trivial task, as no single feature or component is able to distinguish the sets by itself. 

\begin{table}
\begin{tabularx}{\textwidth}{X*{7}{r}}\toprule
          &\multicolumn{3}{c}{Dataset G}      &\multicolumn{3}{c}{Dataset I} \\
\cmidrule(lr){2-4}  \cmidrule(lr){5-7} 
Feature &min	&max	&avg	&min &max		&avg\\
	\midrule
node		          &-0.05	&0.03	&-0.01		&-0.11	&0.08	&-0.02\\
predecessor		      &-0.05	&0.03	&0.01		&-0.43	&0.05	&-0.16\\
\addlinespace
objective		         &-0.14	&0.36	&0.22		&-0.12	&0.36	&0.16\\
consumption\_critical	   &0.03	&0.27	&0.11		&-0.01	&0.26	&0.12\\
tour\_length		           &0.02	&0.37	&0.10		&-0.18	&0.16	&-0.05\\
nunreachable		      &-0.34	&-0.03	&-0.21		&-0.28	&-0.01	&-0.17\\
nvisited		          &-0.34	&-0.06	&-0.21		&-0.28	&-0.03	&-0.17\\
nvisited\_unaltered		  &0.01	&0.38	&0.10		&-0.18	&0.16	&-0.04\\
repeated\_visits		       &-0.04	&0.02	&-0.01		&-0.11	&0.06	&-0.03\\
\addlinespace							
nlabels\_network		&-0.02	&0.20	&0.04		&-0.44	&-0.22	&-0.36\\
nlabels\_open\_network		&0.00	&0.09	&0.03		&-0.17	&0.13	&-0.02\\
nlabels\_closed\_network		&-0.02	&0.17	&0.03		&-0.48	&-0.24	&-0.39\\
nlabels\_dominated\_network		&-0.02	&0.18	&0.04		&-0.44	&-0.24	&-0.37\\
nlabels\_node		&-0.29	&0.19	&-0.11		&-0.37	&-0.15	&-0.28\\
nlabels\_open\_node		&-0.22	&0.03	&-0.10		&-0.22	&0.05	&-0.09\\
nlabels\_closed\_node		&-0.26	&0.12	&-0.08		&-0.39	&-0.19	&-0.29\\

\bottomrule
\end{tabularx}

\caption[Correlation analysis]
{Dataset G and I. Correlation analysis between main features and target class (Pareto-Dominated).} \label{tab:corr}
\end{table}

\begin{figure}
\begin{subfigure}{.49\textwidth}
  \centering
  \includegraphics[width=\linewidth]{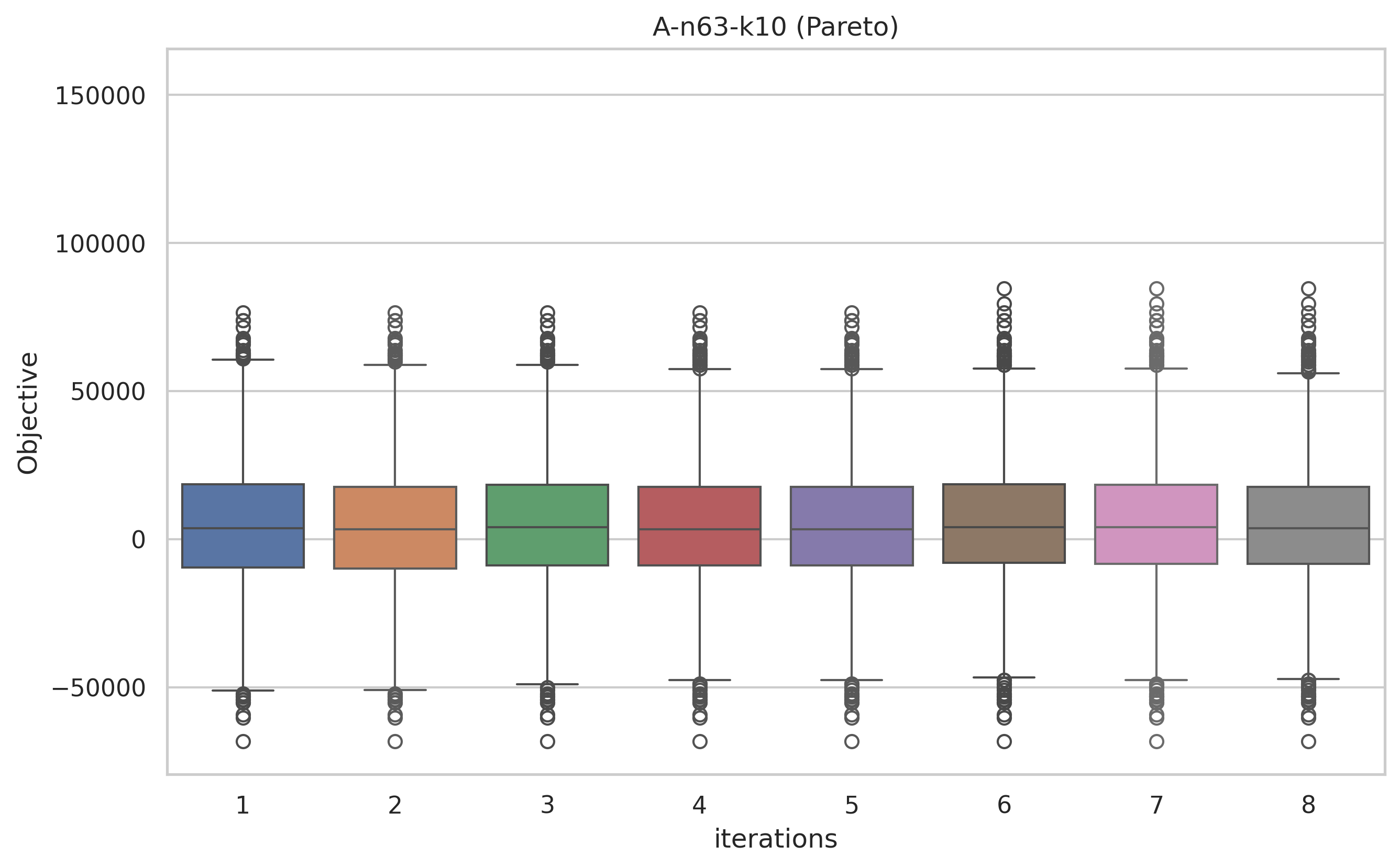}
  \caption{}
  \label{fig:boxobjG1}
\end{subfigure}%
\begin{subfigure}{.49\textwidth}
  \centering
  \includegraphics[width=\linewidth]{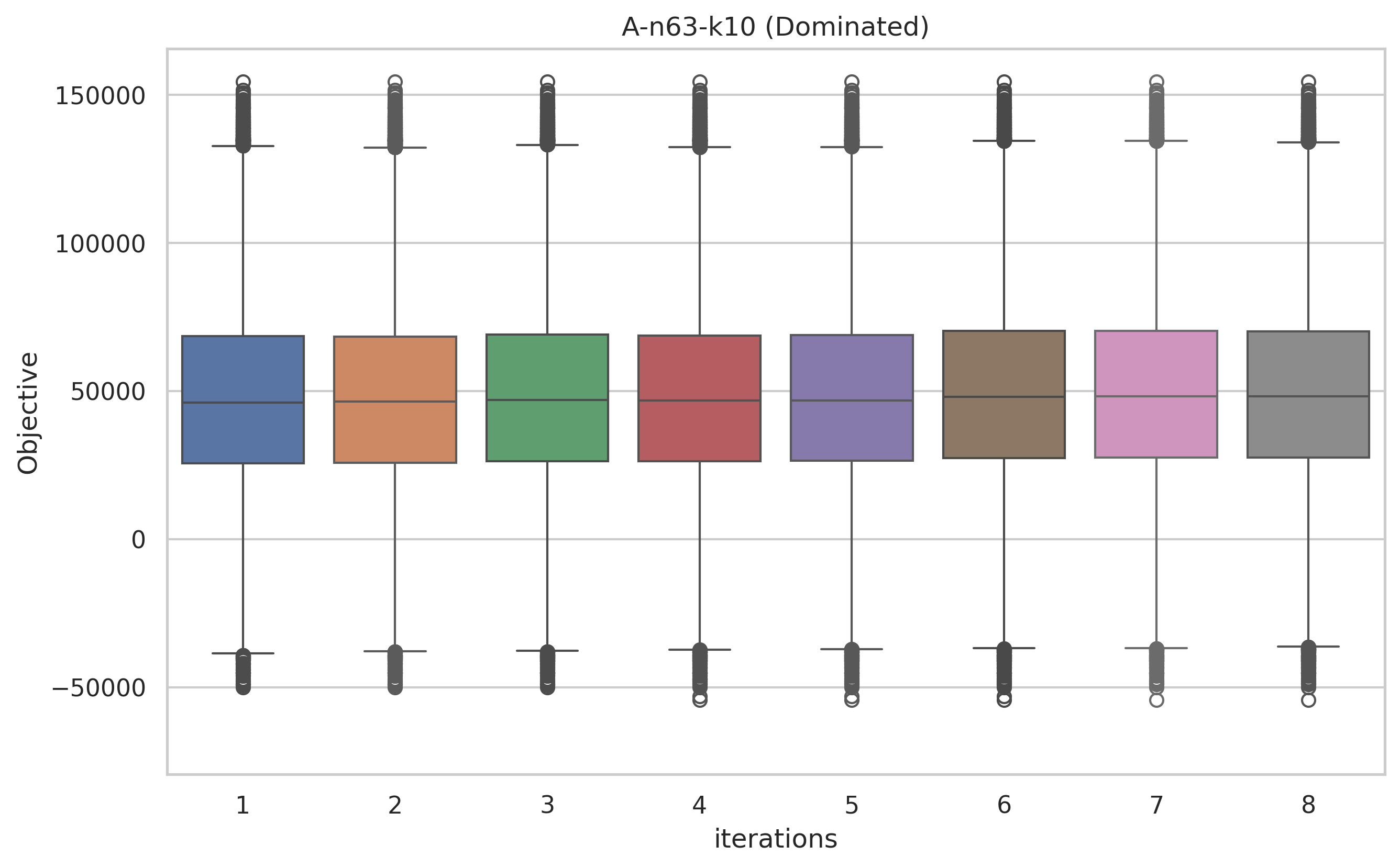}
  \caption{}
  \label{fig:boxobjG2}
\end{subfigure}
\begin{subfigure}{.49\textwidth}
  \centering
  \includegraphics[width=\linewidth]{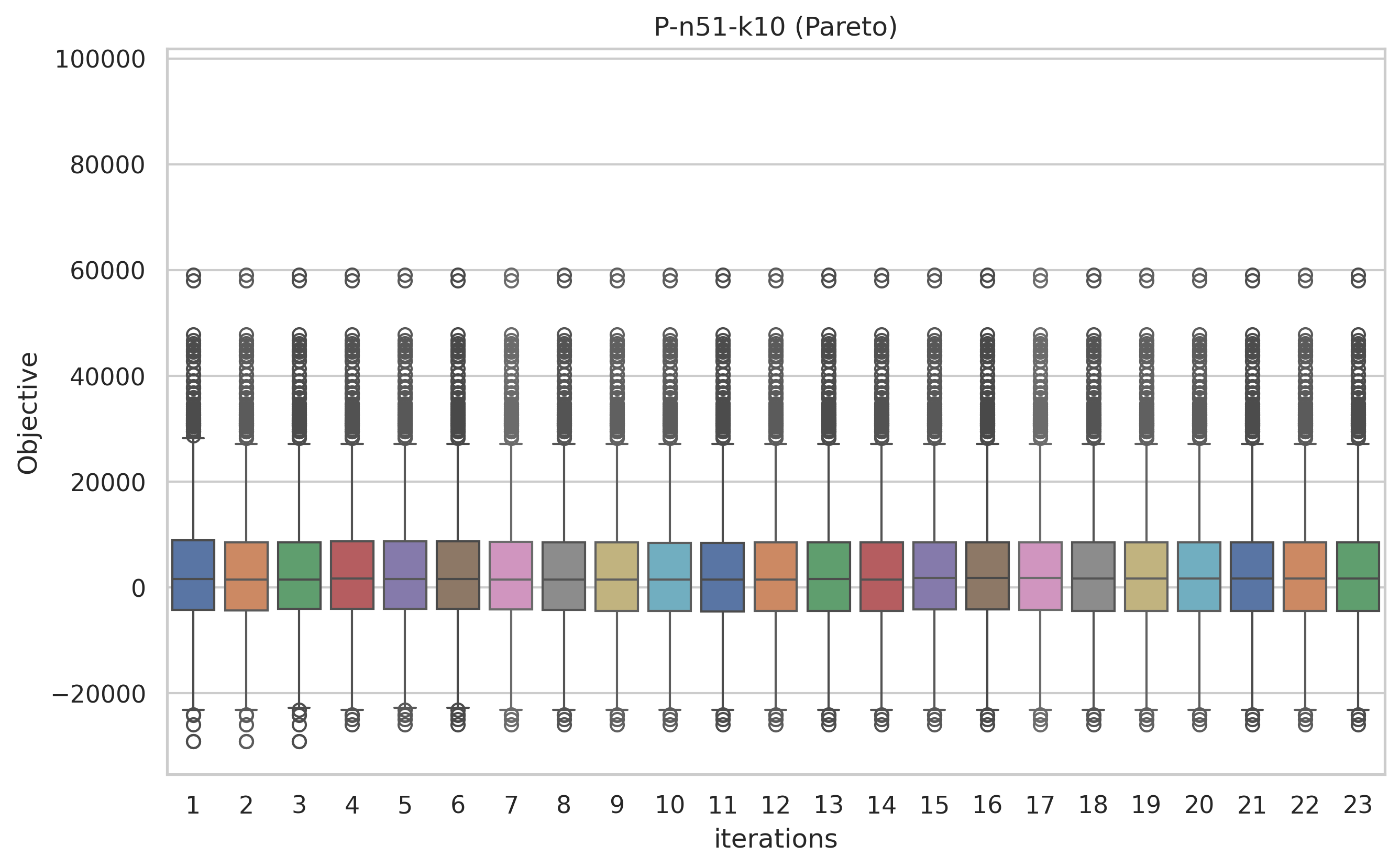}
  \caption{}
  \label{fig:boxobjG3}
\end{subfigure}
\begin{subfigure}{.49\textwidth}
  \centering
  \includegraphics[width=\linewidth]{pictures/dssr/G-A-n63-k10-dom.png}
  \caption{}
  \label{fig:boxobjG4}
\end{subfigure}
\caption{Dataset G. Boxplot of the objective for each iteration of the algorithm. Pareto labels on the left, dominated labels on the right.}
\label{fig: dssr boxplot objective G}
\end{figure}

\begin{figure}
\begin{subfigure}{.49\textwidth}
  \centering
  \includegraphics[width=\linewidth]{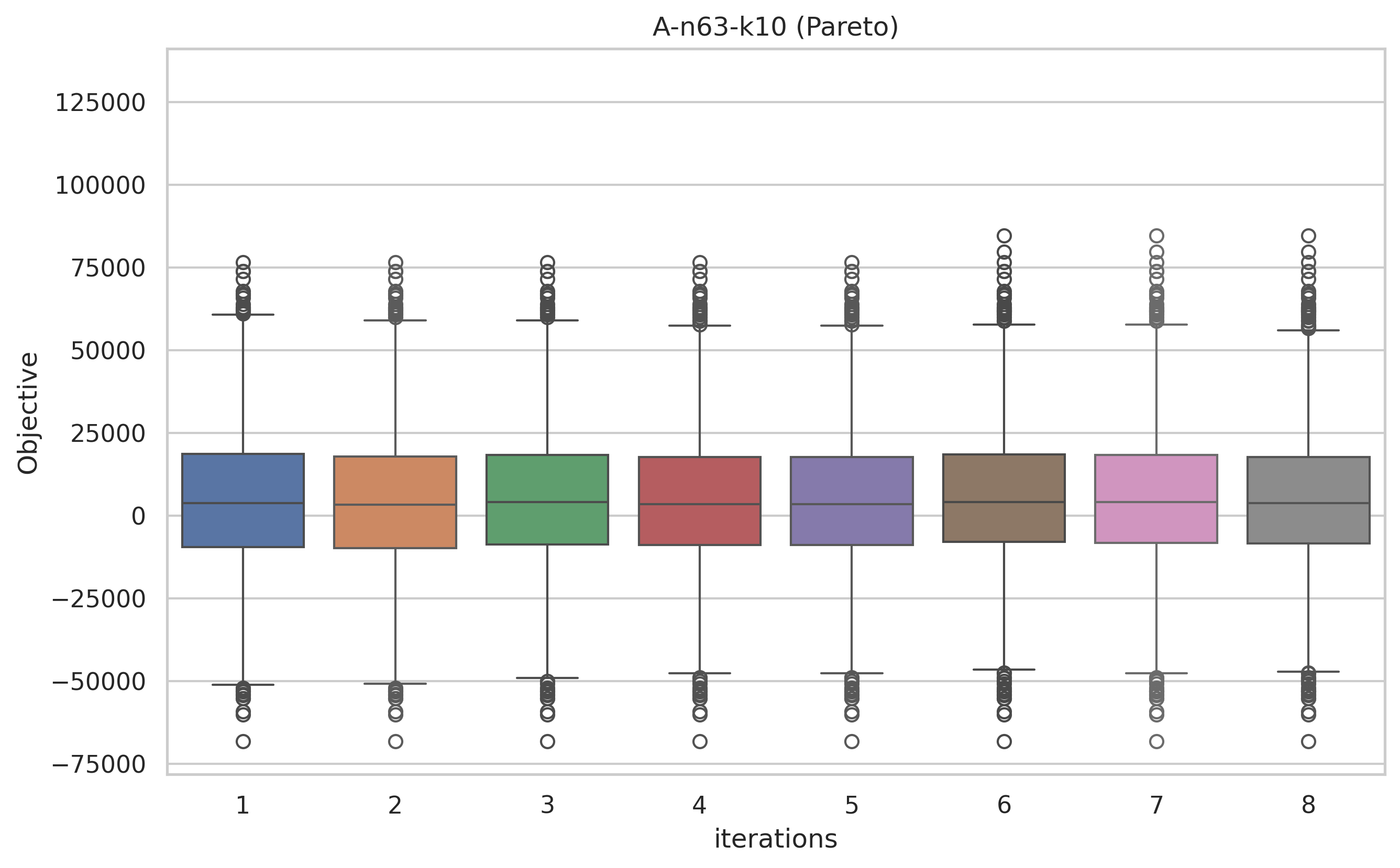}
  \caption{}
  \label{fig:boxobjI1}
\end{subfigure}%
\begin{subfigure}{.49\textwidth}
  \centering
  \includegraphics[width=\linewidth]{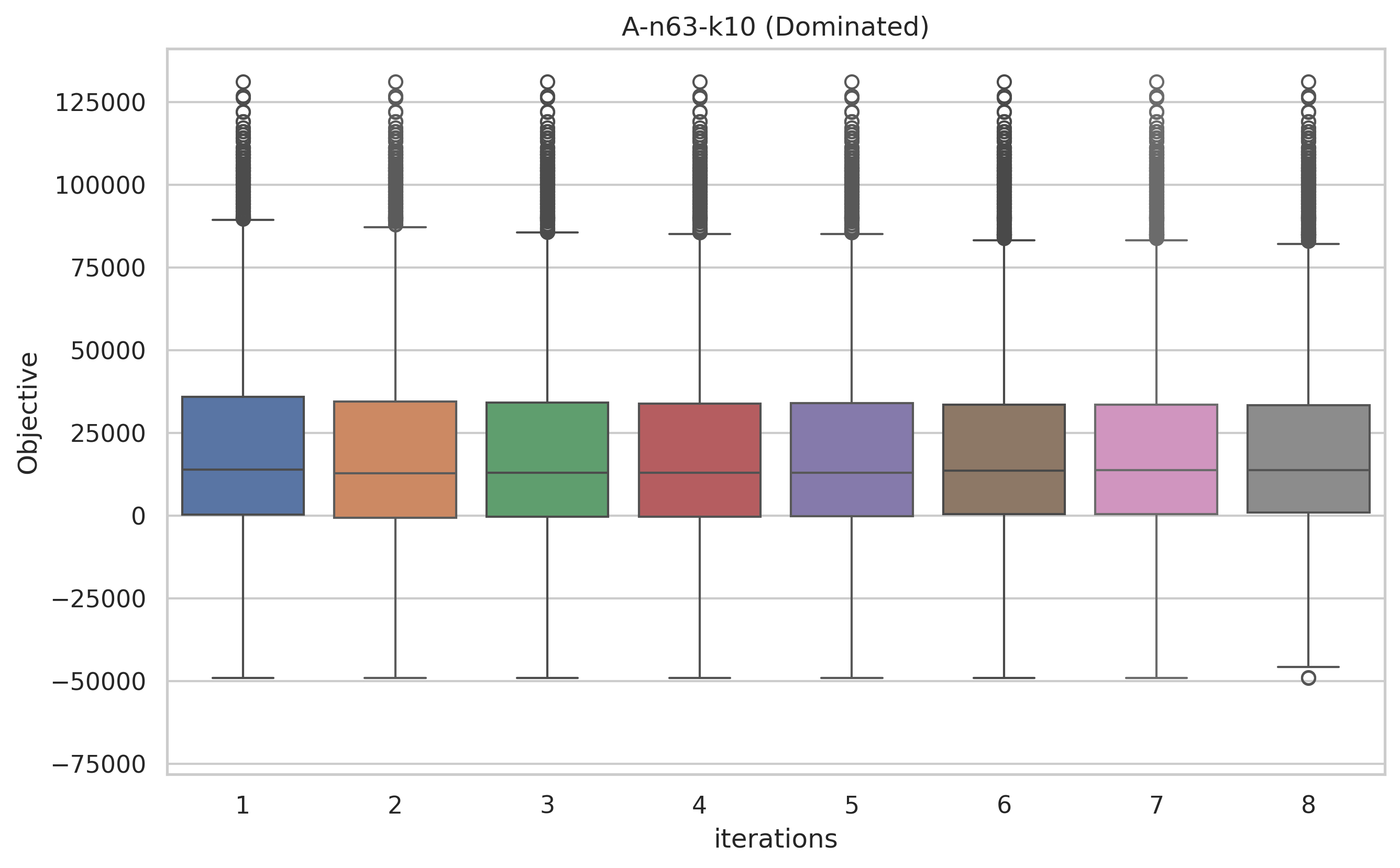}
  \caption{}
  \label{fig:boxobjI2}
\end{subfigure}
\begin{subfigure}{.49\textwidth}
  \centering
  \includegraphics[width=\linewidth]{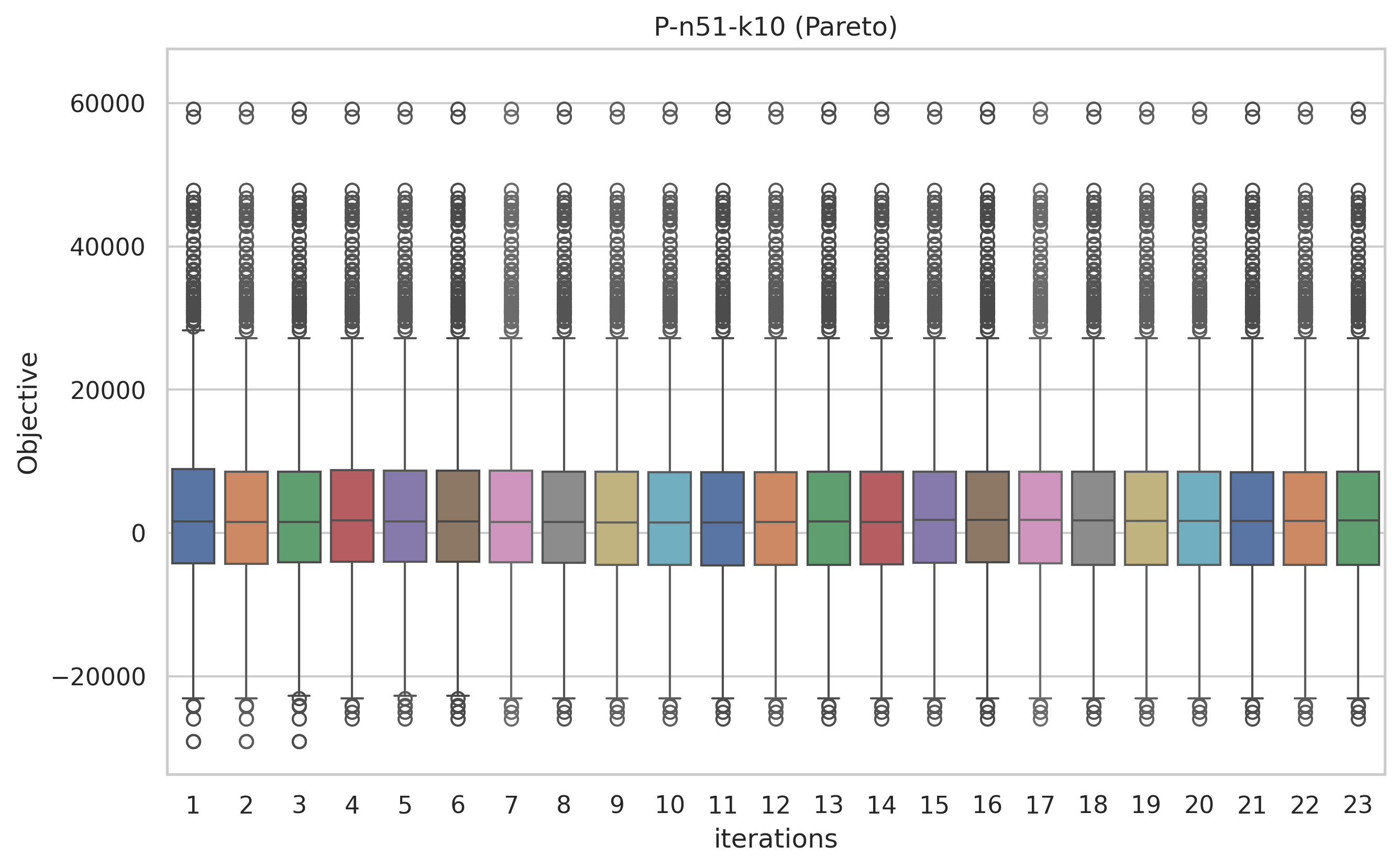}
  \caption{}
  \label{fig:boxobjI3}
\end{subfigure}
\begin{subfigure}{.49\textwidth}
  \centering
  \includegraphics[width=\linewidth]{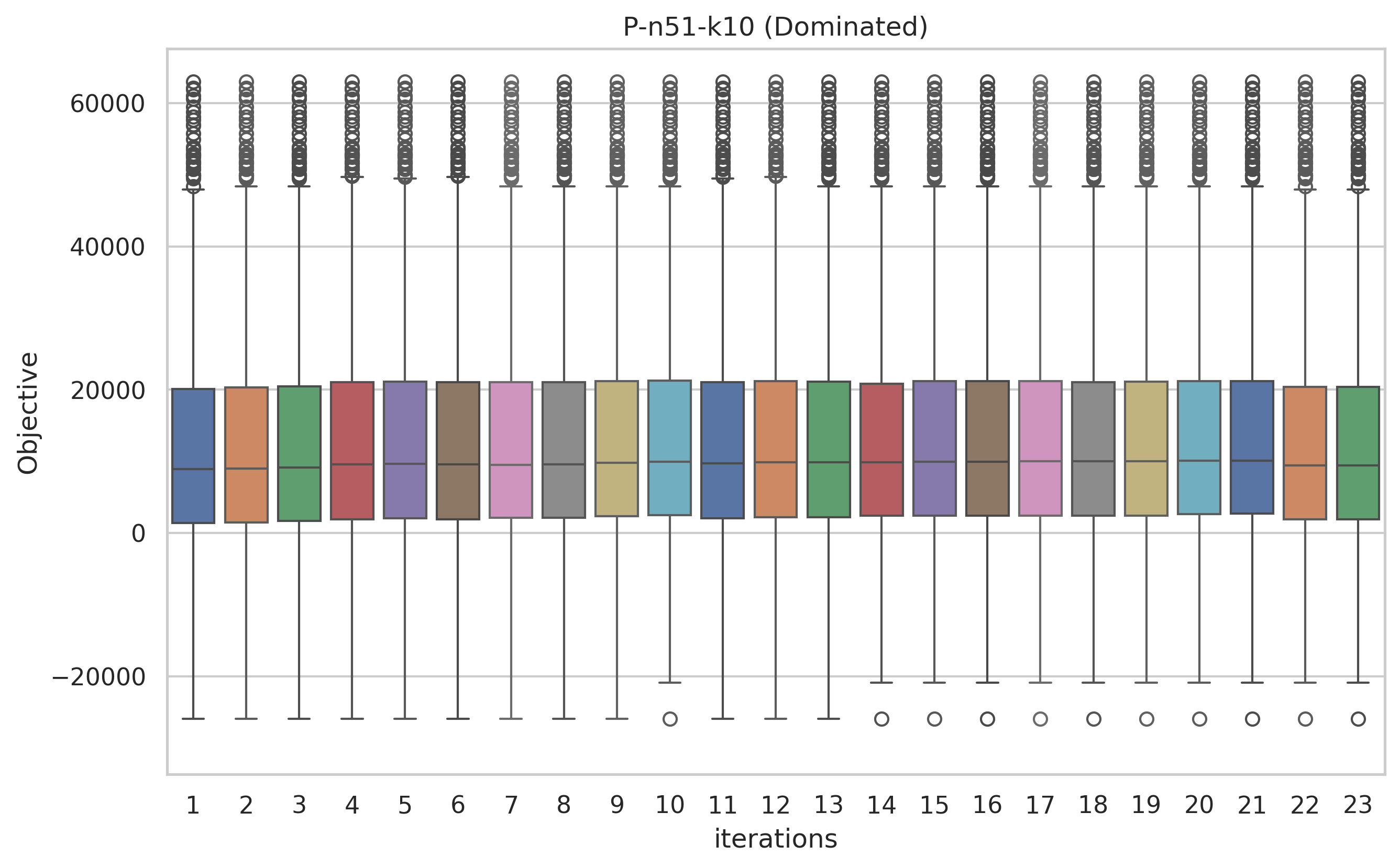}
  \caption{}
  \label{fig:boxobjI4}
\end{subfigure}
\caption{Dataset I. Boxplot of the objective for each iteration of the algorithm. Pareto labels on the left, dominated labels on the right.}
\label{fig: dssr boxplot objective I}
\end{figure}

\paragraph{Analysis 4: DSSR-R iteration analysis}
When employing Decremental State Space Relaxation strategies, algorithms iteratively solve different relaxations of the problem. We follow the intuition that two different subsequent iterations likely share many common labels, and therefore explore, in this analysis, the similarities in the objective and other feature distributions throughout the resolution.
More in detail, we present in Figure \ref{fig: dssr boxplot objective G} a box plot for the objective obtained by each DSSR-R iteration, having solved two selected instances of Dataset G. For each instance, we report the distributions for Pareto labels (left) and for Dominated labels (right). In particular, we observe that the the objective remains similar between iterations: although some small differences are present, this holds true for all problems in SPPRCLIB. However, we note that with other configurations, where the budget of a bidirectional search algorithm gets updated between one iteration and the following one, forward and backward search will likely explore the network differently, therefore producing, overall, a different cost distribution.
We report the same analysis for two selected instances of Dataset I, in Figure \ref{fig: dssr boxplot objective I}. Indeed, the same considerations from before apply, that is, the objective distribution does not differ much, between subsequent iterations. Again, this property remains valid for all instances of SPPRCLIB.
We also repeated the same analysis also for the critical consumption, where we noticed a modest increase in consumption from iteration to iteration: when solving the more relaxed problems, paths with high resource consumption are likely more easily dominated by super optimal labels. %For the latter, we observed that for most problems, the distribution remained largely consistent in most iterations, while in other instances (around 20\% of the total), the tour length distribution changed as DSSR progressed.

Overall, these results suggest that previous iterations likely contain some information that might prove useful in describing the following ones. This might come in the form of already generated labels, or possible upper or lower bound to costs and so on.
Indeed, we detail some of this information during normalization, in the following subsection, and exploit data from previous iterations in some learning setups.

\subsection{Normalization} \label{normalization}

In this section, we explain how we normalized our dataset. Since it contains data from solving various problems, some features, particularly dynamic ones such as the objective value, can vary significantly in scale across instances. Normalization was therefore essential to ensure that features with different scales, also coming from different problem instances, contribute equally during model training. Additionally, normalized data is the requirement for some machine learning models.
Still, normalizing dynamic features like the objective value posed specific challenges due to their strong dependence on the instance context. Indeed, solution-related characteristics can vary widely across instances, influenced by factors such as network topology, cost structures, resource consumption, and more. Consequently, we explored several normalization strategies, carefully weighing their effectiveness and computational impact. 
Therefore, in the following, each instance is normalized using problem-specific data, dynamic information available during label generation, or information derived from the previous DSSR-R iteration, which we accessed in constant time. Indeed, we avoid using information gathered across all instances, such as the lowest bound found or the highest number of nodes, as these may overlook the unique characteristics of new, unseen instances and hinder the generalization to the current ones.

After conducting preliminary experiments, we decided to use problem data to normalize the critical resource consumption. More in detail, we used the upper bound, that is, the overall budget as the normalization factor. For node-related features, such as the number of visited or unreachable nodes, we used the total number of nodes as the reference. 

For normalizing the number of cycles in a path, we used the tour length of the label, instead.
Additionally, we were able to normalize some features related to the network status by using the total number of labels stored in the network at the time of label generation. Specifically, we normalized the number of labels that were yet to be extended (open), already extended (closed), and dominated. At the same time, since we used a queue as a label pool for each node, we also normalized the number of labels that had already been extended (closed) and those yet to be extended (open) with the queue size.

Overall, we found no straightforward solution that used instance specific characteristics for normalizing the objective, the tour length of a path, or the total number of labels stored in the network or queue. Therefore, based on the results of Analysis 4, we explored whether using information from previous DSSR-R iterations of the same instance could be helpful. To this end, we discarded data from the first relaxation of each problem, as no previous relaxation data was available at that point.
%Initially, we attempted to use the best solution from the previous iteration, which serves as a lower bound for the current one, to normalize the objective. However, this approach proved unsuitable, as there is no certainty regarding how much the overall bound will change between iterations. As a result, Pareto labels present in both iterations, particularly the one generated at the beginning of the algorithm, would have different objective scores. %Additionally, this bound is derived from combining a forward and a backward label, while the majority of labels in the set correspond to partial paths only. 
After some experiments, we found that employing a min-max normalization by using the objective of the partial paths obtained in the previous iteration worked best for our goal. More in detail, for a label $l$ in relaxation iteration $i$, we used the objective from the label set $L$ of iteration $i - 1$  to normalize it.
\begin{equation}
obj_{l, norm} = \frac{obj_l - \min_{L_{i-1}} \{ obj \}}{\max_{L_{i-1}} \{ obj \} - \min_{L_{i-1}} \{ obj \}}
\end{equation}

\begin{figure}
\begin{subfigure}{.49\textwidth}
  \centering
  \includegraphics[width=\linewidth]{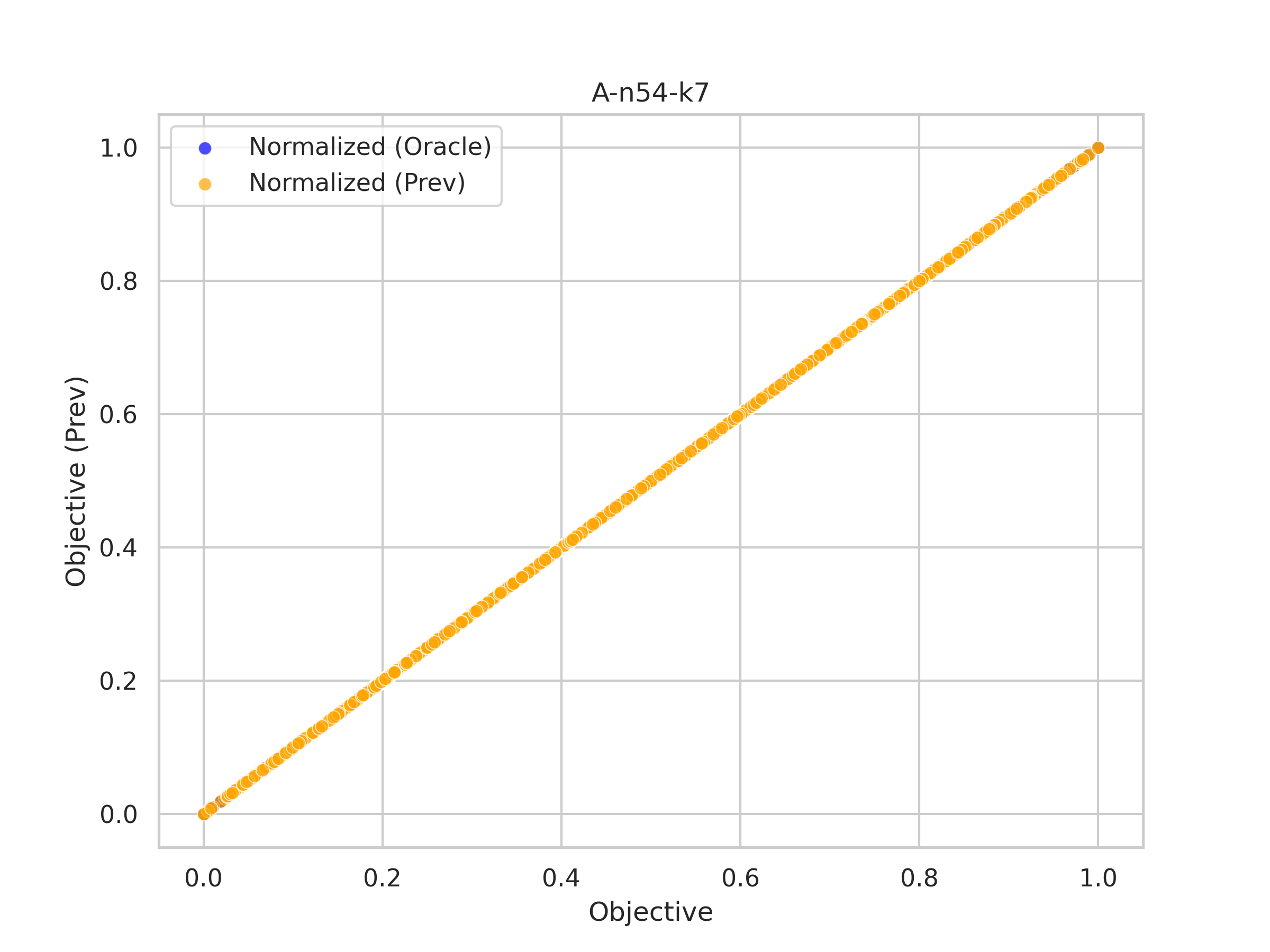}
  \caption{}
  \label{fig:normGA}
\end{subfigure}%
\begin{subfigure}{.49\textwidth}
  \centering
  \includegraphics[width=\linewidth]{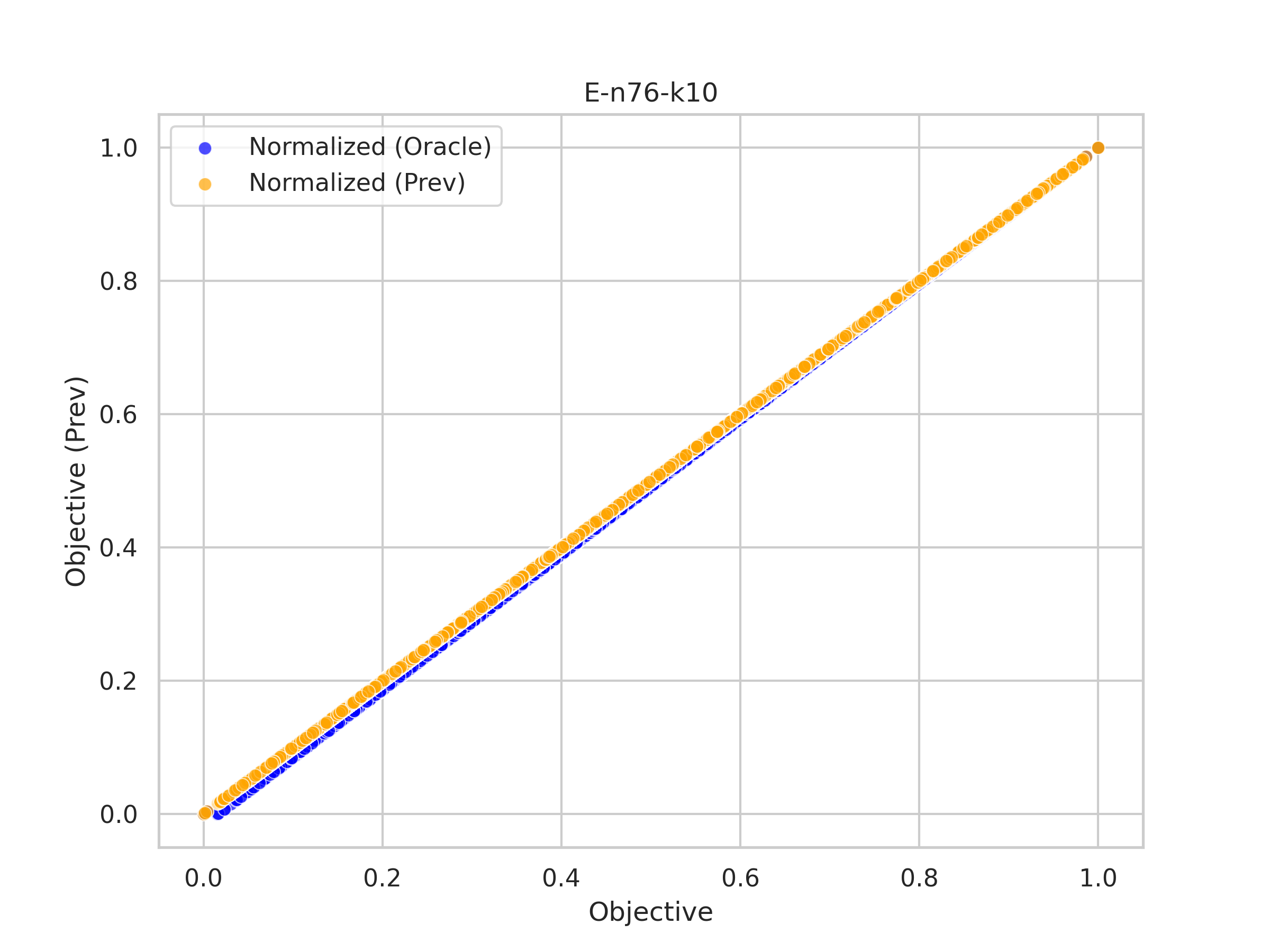}
  \caption{}
  \label{fig:normGE}
\end{subfigure}
\caption{Dataset G. Plotting objective normalized with previous data (Prev) against objective normalized with current data (Oracle).}
\label{fig: normalizationG}
\end{figure}

\begin{figure}
\begin{subfigure}{.49\textwidth}
  \centering
  \includegraphics[width=\linewidth]{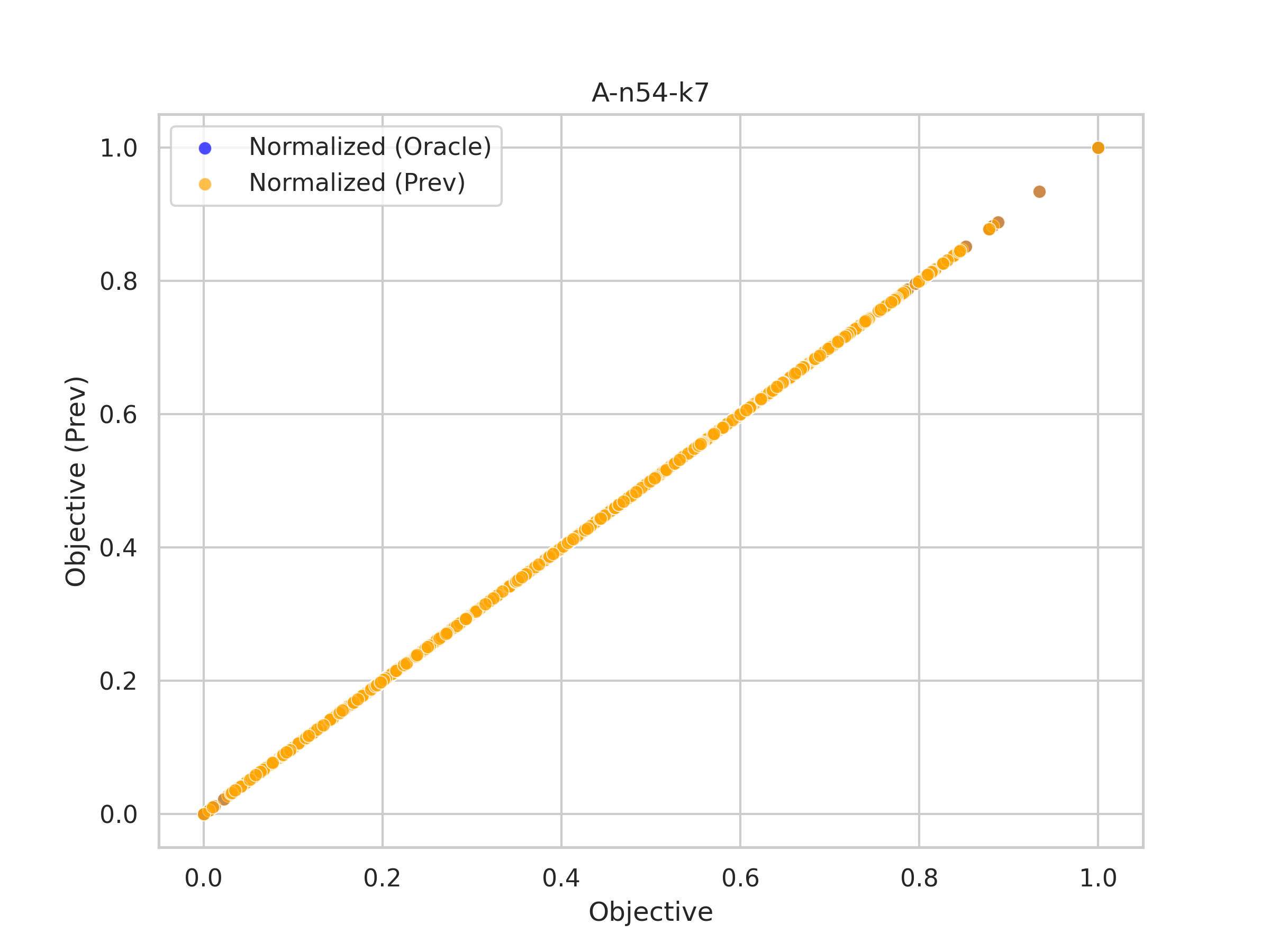}
  \caption{}
  \label{fig:normIA}
\end{subfigure}%
\begin{subfigure}{.49\textwidth}
  \centering
  \includegraphics[width=\linewidth]{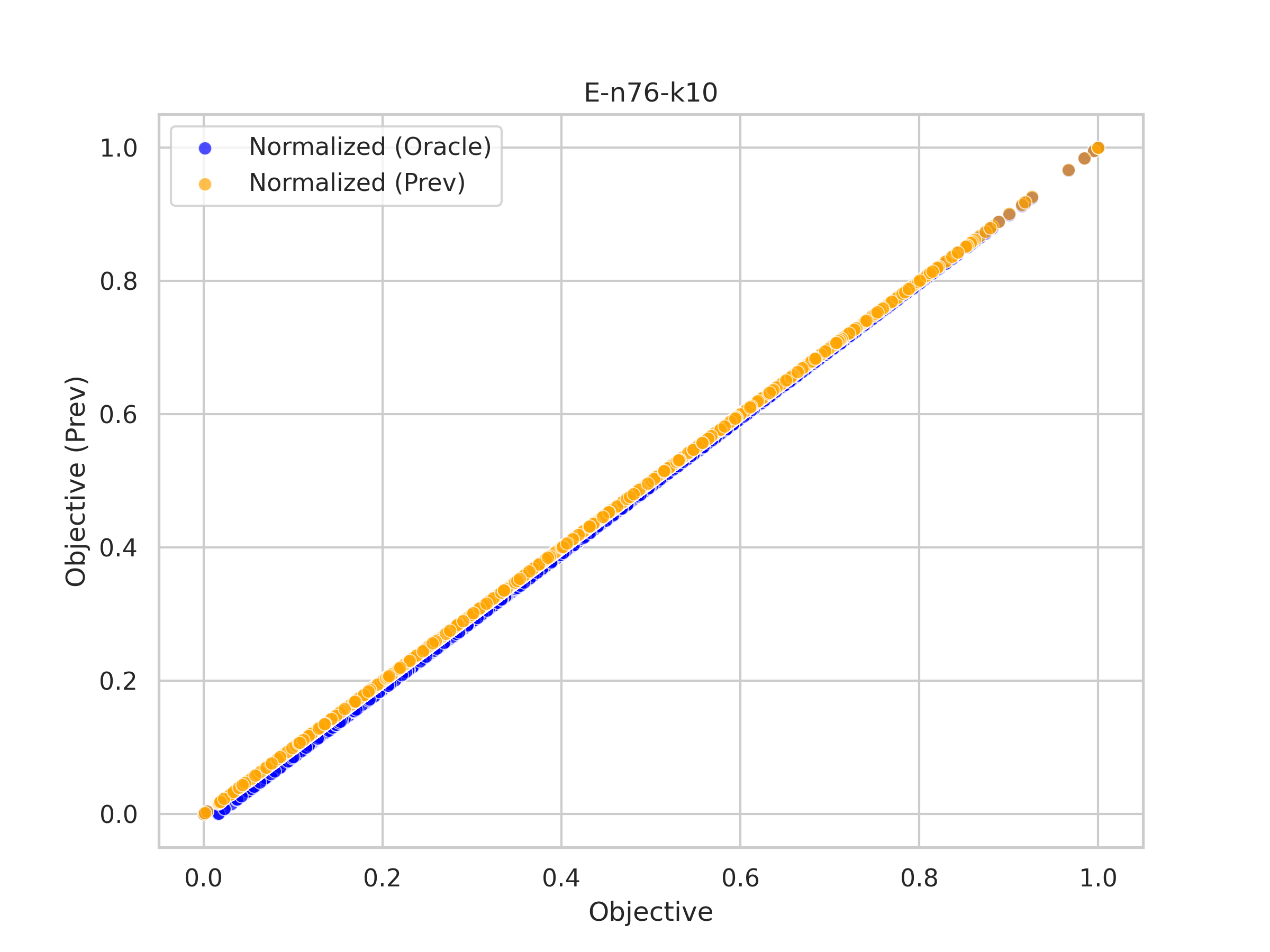}
  \caption{}
  \label{fig:normIE}
\end{subfigure}
\caption{Dataset I. Plotting objective normalized with previous data (Prev) against objective normalized with current data (Oracle).}
\label{fig: normalizationI}
\end{figure}

For validation purposes we also performed a min-max normalization by using the full label set of the current iteration. We call this normalization Oracle, as it can be only be performed a posteriori, once the relaxation has been fully solved.
We plot our normalized objective (Prev, in orange) against the Oracle normalization (in blue) for the last iteration of selected instances in Figure \ref{fig: normalizationG} and Figure \ref{fig: normalizationI} for Dataset G and I, respectively.

Overall, we observed no differences in normalization scores for all instances but two. That is, the majority of instances shows a scatter plot similar to \ref{fig:normGA} and \ref{fig:normIA}. For the remaining two, the error of normalized scores was less than 0.02 (see Figures \ref{fig:normGE}, \ref{fig:normIE}).

We found this approach to normalization to be suitable. Alternative strategies, including the z-score, were also considered and showed similar performance. However, we note some drawbacks: if the search budget is updated at the end of a DSSR-R iteration, the objective distribution in the next iteration will likely change, which in turn may affect the normalization process. This issue, however, can be avoided by selecting a fixed budget at initialization, possibly after some pre-processing, and by preventing dynamic updates. 
Overall, this strategy proves particularly effective because DSSR-R updates the bitsets of only a subset of the network nodes at each iteration. In principle, this means that each relaxation step is similar to the previous ones. Preliminary experiments confirm that this normalization strategy remains effective even when applied to standard DSSR methods, although the resulting accuracy is lower.

We identified no consistent technique to normalize the number of labels and the tour length. For the number of labels, we used the number obtained after solving the previous relaxation. However, since DSSR-R solves progressively less relaxed problems, the number of labels tends to increase from one iteration to the next one: as a result, achieving a 0 to 1 normalization was not possible. For the tour length, we used the length of the longest path found in the previous optimization round as the normalization factor. While one could consider using the overall number of nodes as an alternative, this could pose issues in large networks with short paths, where such a normalization would compress the scores, resulting in values close to 0.

Finally, we also designed an additional feature, called efficiency, as a possible proxy of the quality of the label. For a given label $l$ we compute it as follows:
\begin{equation}
eff_l = 1 - obj_{l, norm} * consumption\_critical_{l,norm}
\end{equation}
That is, the efficiency score is high if both the objective and the consumption critical are low. We included this feature in subsequent analyses that use normalized features, still, we note that we found no strong direct correlation with the label type.

\section{Learning to discriminate labels} \label{sec:learning}
Following the results of the previous analyses, we turned our attention to a more complex question: whether combinations of features exist that can distinguish between Pareto and Dominated Labels, and, if so, to understand which setups make these classifiers work effectively.
Therefore, we explored various policies to train supervised models for both Datasets G and I. 

We propose settings that incorporate inter-instance models, as well as those that operate strictly within individual instances.

In the first case, we consider a scenario where a set of problems has already been solved, their labels have been collected, and one or several models have been trained offline using different strategies. These models are then tested on the labels of new, unseen problems. We report this case in subsection \ref{offline-learning}.

In contrast, the second setup is presented in \ref{online-learning} and it involves facing a completely new problem, when no pre-existing models are available. In this case, we simulate having access to only a subset of the available data from the current problem to train a model online, on the fly, and then use it to distinguish upcoming promising labels.

\subsection{Training strategy}
In this subsection, we describe the policy that we used in \ref{offline-learning} and \ref{online-learning} to train our models.
First, we always trained using data from a single iteration of DSSR, typically the latest available at run time. In what follows, we describe the setup for training models on both single instances and aggregated groups of instances.

Given that the class distribution can be very imbalanced, particularly in Dataset G, we applied sampling to ensure a 1:1 ratio between Pareto and Dominated labels for each instance. Specifically, we sampled approximately 80\% of the number of labels in the minority class and matched that count with an equal number of labels from the majority class. The remaining labels were reserved for validation purposes. Without proportional sampling, models trained on Dataset G would tend to be biased toward predicting Dominated labels.

When training on data coming from multiple instances, we took into account that the total number of generated labels can vary significantly across them. To prevent a single instance from dominating the training data, we limited the size of the sample, for each instance, by using the median number of labels from the minority class across the group.

After preliminary experiments, we chose to train XGBoost (version 2.1.1) classifiers using the default settings. The {\tt dominated} feature (0 -- Pareto, 1 -- Dominated) served as the label for the classifier, while the features marked with an asterisk in Section~\ref{appendix:features} of the appendix were used as input.  This model proved suitable for our needs, as it accommodates both normalized and non-normalized data. While preliminary tuning did not yield significant gains, a more systematic exploration of hyperparameter optimization could potentially warrant further investigations.
Overall, we found this setup to be the best performing one in our early evaluations.

We discuss additional details regarding training also in the following subsections, while full validation results are reported in the Appendix, in \ref{appendix:validation}.

\subsection{Offline learning} \label{offline-learning}
We focus on the first ambitious task, that is understanding if, given the labels generated by solving an instance or a collection of instances, it is possible to discriminate among Pareto and Dominated labels of a previously unknown problem. To this end, we trained our models by aggregating instances in accordance to two main strategies: one-vs-all and all-but-one. Additional instance aggregation policies are summarized in the Appendix, in Section \ref{appendix:aggregation}.  In the following, all training and testing was performed on normalized data from the last relaxation of the involved instances. 

\paragraph{One-vs-All strategy}
First, we trained a single model for each instance, and used it to predict the label type of all the remaining ones, in a standard one-vs-all fashion.

We present the test results in Table \ref{tab:offline-single} for both datasets, reporting prediction accuracy for both label types, precision, and F1-score. We remark that in our setup, Recall corresponds to the accuracy on Dominated labels. In all these metrics, higher values are preferable. Furthermore, all scores are averaged based on the training model’s class. For example, class A represents the performance of models trained on that class in predicting the label type of all other instances. This includes instances from the same class.

We also reported the average absolute difference between the two accuracies (Diff.), computed per instance and then averaged across all instances in the group. This metric captures the prediction bias of the models toward one of the two classes: a score of 0 means the model classifies Pareto and Dominated labels with equal accuracy on that instance, while a score of 1 indicates the model predicts only one of the two classes for that instance.

We observe that prediction is possible on Dataset G, achieving a balanced accuracy of 78\%. However, while all models generally perform well in identifying the majority class (Dominated labels), their performance on the minority class (Pareto-labels) is notably lower. Specifically, the Diff. column indicates a score difference of about 0.3, strongly highlighting the class imbalance in predictions. Overall, models trained on problems from class A and B yield the most balanced results, whereas those trained on class E and M were significantly more biased toward dominated labels.

Interestingly, we also found that classification ability was not symmetrical. Let $\delta$ and $\gamma$ be two different instances, with models trained on their records denoted as $M\delta$ and $M\gamma$, respectively.  If $M\delta$ performs well in predicting $\gamma$, there is no guarantee that $M\gamma$ will perform well in predicting $\delta$.

On the other hand, tests on Dataset I reveal that distinguishing a Pareto-label from a Dominated label, in the label pool, is significantly more challenging. No model achieves satisfactory performance, and predictions exhibit substantial differences in class scores. This difficulty arises because the labels in the label pools are intrinsically more similar to each other. Additionally, Dataset I contains fewer records, and smaller instances likely lack sufficient information for training to generalize effectively to other instances.

Ultimately, while this strategy appears feasible when enough data is available, individual instances alone do not guarantee accurate generalization. Consequently, we explored alternative strategies.

\paragraph{All-but-one strategy}
Therefore, in our next experiment, we trained an XGBoost model using data from all instances except one, which was reserved for testing. In other words, for each test instance, we trained a model incorporating information from all other instances. This setup simulates a scenario where a model has already been trained, and a completely new, unseen problem, needs to be analyzed.

We present the results in Table \ref{tab:offline-all} for both datasets. Specifically, we report the average class-wise accuracies, the average absolute difference between them (Diff.), precision, and F1-score. In this case, average scores are aggregated based on the class of the testing instance. Therefore, every row describes how well the models are able to predict a single instance belonging to that particular class.

Overall, for offline learning, this approach yielded the best results. In particular, the results on Dataset G show higher accuracy (85.5\% balanced accuracy) and lower Diff. values compared to other strategies. 
Indeed, this policy enables prediction across all instance classes in Dataset G: on average, no class shows strong bias in its scores.
 
However, the same does not fully hold for Dataset I. While this approach achieved the best results among offline strategies (around 75\% balanced accuracy and much lower Diff. values) this level of accuracy, combined with the still non-negligible Diff. values (albeit significantly lower than those of the one-vs-all strategy), remains substantial enough to impact the quality of predictions.

\paragraph{Perspectives} Overall, these findings suggest that distinguishing label types (Dataset G) in previously unseen problems can be performed with adequate accuracy while classification of labels stored in the pool (Dataset I), is feasible but inconsistent. Indeed, training more accurate predictors for both datasets remains an open research challenge. This may require designing and integrating more dynamic features, which could involve pre-processing steps or, more generally, non-constant computational times: as a result, there may be a trade-off between accuracy and computational efficiency.

Finally, we observe that using a fixed 1:1 sampling ratio between label types and a single training configuration may not be suitable for all test instances. A promising direction for improvement would be to introduce pre-processing techniques and cluster similar instances together. In our experiments, we also attempted clustering based on class, that is, grouping problems generated in a similar manner, but this did not lead to the best predictive performance. Developing customized models tailored to specific groups of similar problems, characterized by more than just their instance generation strategy, could prove to be a more effective approach.

\begin{table}
\begin{tabularx}{\textwidth}{X*{8}{r}}\toprule
        &    Training          &  Test      &\multicolumn{3}{c}{Accuracy} & & &\\
\cmidrule(lr){4-6} 
Dataset &Model	& Inst.	&Pareto	&Dominated &Diff.		&Precision	&F1\\
	\midrule
    \multirow{4}{*}{G}
&A	&400	&0.78	&0.83	&0.26	&0.99	&0.89\\
&B	&280	&0.81	&0.79	&0.26	&0.99	&0.86\\
&E	&240	&0.59	&0.90	&0.40	&0.98	&0.93\\
&M	&120	&0.55	&0.92	&0.41	&0.98	&0.94\\
&P	&600	&0.71	&0.80	&0.36	&0.99	&0.86\\

\addlinespace
Overall &&	1640	&0.72	&0.83	&0.33	&0.99	&0.88\\
\hline
\addlinespace
\multirow{5}{*}{I}
&A	&400	&0.69	&0.63	&0.41	&0.71	&0.60\\
&B	&280	&0.78	&0.55	&0.43	&0.75	&0.57\\
&E	&240	&0.56	&0.76	&0.45	&0.66	&0.66\\
&M	&120	&0.39	&0.85	&0.55	&0.59	&0.67\\
&P	&600	&0.75	&0.55	&0.50	&0.76	&0.55\\

\addlinespace
Overall &&	1640	&0.69	&0.62	&0.46	&0.72	&0.59\\
\bottomrule
\end{tabularx}

\caption[Inter-class]
{Dataset G and I. Classification results, offline learning. Classifiers were trained on a single instance and tested against the others.} \label{tab:offline-single}
\end{table}

\begin{table}[h!]
\small

\begin{tabularx}{\textwidth}{X*{8}{r}}\toprule
        &      \multicolumn{2}{c}{Test}                &\multicolumn{3}{c}{Accuracy} & & &\\
\cmidrule(lr){2-3} \cmidrule(lr){4-6} 
Dataset &Class	&Inst.	&Pareto	&Dominated &Diff.		&Precision	&F1\\
	\midrule
    \multirow{4}{*}{G}
&A	&10	&0.86	&0.86	&0.17	&0.99	&0.92\\
&B	&7	&0.80	&0.83	&0.19	&0.99	&0.89\\
&E	&6	&0.97	&0.84	&0.13	&1.00	&0.91\\
&M	&3	&0.91	&0.89	&0.06	&1.00	&0.94\\
&P	&15	&0.84	&0.85	&0.18	&0.99	&0.91\\

\addlinespace
Overall&	&41	&0.86	&0.85	&0.16	&0.99	&0.91\\
\hline
\addlinespace
\multirow{5}{*}{I}
&A	&10	&0.78	&0.72	&0.18	&0.77	&0.73\\
&B	&7	&0.71	&0.78	&0.14	&0.75	&0.76\\
&E	&6	&0.86	&0.63	&0.25	&0.84	&0.69\\
&M	&3	&0.8	&0.69	&0.26	&0.88	&0.77\\
&P	&15	&0.73	&0.8	&0.17	&0.7	&0.73\\
\addlinespace
Overall &	&41	&0.76	&0.74	&0.18	&0.76	&0.73\\
\bottomrule
\end{tabularx}

\caption[Inter-class]
{Dataset G and I. Classification results, offline learning. For each experiment, the classifier was trained on all instances except one, that was used for testing.} \label{tab:offline-all}
\end{table}

\subsection{Online learning} \label{online-learning}
Next, we tried to determine whether Pareto labels could be distinguished from Dominated ones within a given instance and a given relaxation, by using historical labels already obtained for that instance during optimization. Indeed, one possibility would be to use labels produced at runtime, during the current relaxation, to make an online prediction of the subsequent labels. While this is theoretically possible, it introduces several algorithmic challenges, such as managing additional computational overhead and understanding how many labels are necessary to train accurate models. Furthermore, there is no guarantee that a set of initial labels will contain enough data to make accurate predictions for the remaining ones.
Indeed, preliminary experiments showed that the models performed well only when most of the states had already been generated. In this regard, we note that our models also rely on features that describe the state of the network and data structures at the moment a label is created or inserted. In scenarios where only partial information is available, these features lose much of their value, as they do not capture the full evolution of the search process. In other words, the model lacks access to the complete history of the network exploration.

Given the results observed in Analysis 4, we therefore opted for a more operationally sensible scenario: assuming that all labels produced by solving relaxation $i$ are already available, we aim to make informed decisions when addressing the subsequent relaxation $i+1$. Specifically, we use the data from relaxation $i$ to train a supervised model and classify labels from iteration $i+1$. 

Indeed, we trained an XGBoost model using labels from the penultimate relaxation solved with DSSR-R. We then tested this model on all the completely unseen labels coming from the last relaxation, that is, the iteration of the algorithm that provides an elementary solution. It is worth noting that the last iteration of DSSR-R solves the least relaxed problem, making it the most challenging and time consuming of the algorithm. 
Since all data originated from the same instance, we trained the model using non-normalized records and, given the limited number of labels available in a single iteration, training time was relatively short, averaging, over all instances, 0.79 seconds for Dataset G and 0.11 seconds for Dataset I, making the approach suitable for this online scenario.

We report test results in Table \ref{tab:online} for both Dataset I and G. We grouped all test instances by their respective class and, for each group, we reported the number of instances (Inst.) along with the average accuracy in classifying Pareto and Dominated labels. Additionally, we included the Diff. score, Precision and F1-score.

\begin{table}[h!]
\small

\begin{tabularx}{\textwidth}{X*{8}{r}}\toprule
        &       \multicolumn{2}{c}{Test}               &\multicolumn{3}{c}{Accuracy} & & &\\
\cmidrule(lr){2-3} \cmidrule(lr){4-6} 
Dataset &Class	&Inst.	&Pareto	&Dominated &Diff.		&Precision	&F1\\
	\midrule
    \multirow{4}{*}{G}
&A	&10	&0.98	&0.93	&0.06	&1.00	&0.96\\
&B	&7	&0.98	&0.90	&0.08	&1.00	&0.95\\
&E	&6	&0.97	&0.94	&0.06	&1.00	&0.97\\
&M	&3	&0.98	&0.93	&0.05	&1.00	&0.97\\
&P	&15	&0.97	&0.92	&0.05	&1.00	&0.96\\

\addlinespace
Overall&	&41	&0.98	&0.92	&0.06	&1.00	&0.96\\
\hline
\addlinespace
\multirow{5}{*}{I}
&A	&10	&0.90	&0.87	&0.04	&0.90	&0.89\\
&B	&7	&0.89	&0.83	&0.06	&0.89	&0.86\\
&E	&6	&0.91	&0.88	&0.03	&0.91	&0.90\\
&M	&3	&0.88	&0.81	&0.07	&0.93	&0.86\\
&P	&15	&0.90	&0.89	&0.06	&0.88	&0.88\\
\addlinespace
Overall &	&41	&0.90	&0.87	&0.05	&0.89	&0.88\\
\bottomrule
\end{tabularx}

\caption[Intra-class, Inter-iteration]
{Dataset G and I. Classification results, online learning. Model trained on iteration $i$ data, tested on iteration $i+1$ records.} \label{tab:online}
\end{table}

In dataset G, performance is consistently high across all classes. Both Pareto and Dominated labels are classified with high accuracy (above 0.90 in all cases), and the low Diff. values (that range from 0.05 to 0.08) indicate a balanced treatment of both classes. Precision and F1-scores are very strong: overall, this suggests that the classifier performs reliably across all classes in this dataset, with low bias toward one label over the other.

Dataset I, while showing lower overall accuracy, demonstrates strong predictive capability. Indeed, instances are classified with high accuracy (on average, above 0.88 balanced accuracy), and although Dominated accuracy is slightly lower, particularly in class M, other classes (like A, E and P) show very balanced performance. 

The average Diff. is still very low (0.05), reflecting some minor bias toward the Pareto class. Nevertheless, the model achieves good Precision (0.89) and F1-score (0.88), indicating that while predictions are less accurate than in Dataset G, they remain reliable and practically useful.

We note that preliminary experiments were also conducted on datasets generated with DSSR (instead of DSSR-R) and they yielded similar results, but with a moderate drop in performance for Dataset I. While the proposed strategy is applicable to any DSSR setup, it is especially effective with DSSR-R, where the similarity between consecutive iterations is higher due to modifications being restricted only to the sets of the nodes involved in cycles.

In summary, the classifiers performed very well on both datasets, demonstrating that it is possible to learn online to distinguish between Pareto and Dominated labels when historical label data from previous relaxations is available. We expect that even better results could be achieved through instance specific tuning during training.

\section{Conclusions and future research} \label{sec:conclusions}

In this paper, we proposed a first exploration of data driven approaches for combinatorial optimization problems, that manage information obtained by dynamic programming states during iterative relaxations optimization. In particular, we focused on shortest path problems with resource constraints and labeling algorithms, and addressed the possibility of discriminating label types, that is Pareto-optimal states from suboptimal ones, looking only at their static properties.

To this end, we designed computationally efficient features and stored every label generated from solving all SPPRCLIB instances with a state-of-the-art solver, PathWyse, thus generating two datesets. The former, Dataset G, includes all the generated labels (P,GD,GI), while the latter, Dataset I, collects only labels that have been inserted in data structures (P, GI). That is, suboptimal states, that are discarded at generation by the algorithm, are not saved in Dataset I.

We then conducted an exploratory analysis of the main features, confirming that labeling algorithms prune suboptimal states efficiently. Although we found no single strongly correlated feature to distinguish between the states, the distribution of the objective value alone can provide some guidance for dataset G, indicating that the states may still be meaningfully separated. Notably, objective distribution also remains largely consistent across DSSR-R iterations for both datasets, suggesting that previously solved relaxations may contain useful information for subsequent ones. We found similar properties with resource consumption. Therefore, we devised a feature normalization strategy based on these insights.

Finally, we trained, for each dataset and each instance, an XGBoost model for binary classification of previously unseen labels, using a 1:1 data split, ensuring an equal number of Pareto and dominated labels. 

First, we simulated offline learning setups, where models were trained using labels from previously solved instances, applying various aggregation strategies. These models were then evaluated on unseen problems and labels. For Dataset G, the predictions showed good accuracy under an all-vs-one policy, though with some bias toward a single class. In contrast, the models struggled to accurately generalize and classify Dataset I labels, highlighting that building more effective models for Dataset I remains an open challenge. Indeed, distinguishing between the subset of Pareto and dominated labels inserted into label pools is significantly more difficult, as these labels are inherently more similar.

Next, we examined scenarios simulating online learning, where models are trained using data from the previous relaxation step. These settings yielded strong predictive performance, with approximately 90\% balanced accuracy for both Dataset G and Dataset I. This confirms that the model can effectively learn from the label sets produced in the preceding DSSR-R iteration.

Still, determining the extent to which these results stem from the models generalization capabilities versus their ability to recognize patterns from previously seen records requires further investigation. Nonetheless, the effectiveness of these methods is well-supported by the results.

Overall, our work demonstrates that, on SPPRCLIB, predicting Pareto labels from Dominated ones is feasible, both in an online and offline setting, achieving good accuracy even on previously unseen problems. This research provides a first foundation for analyzing dynamic programming states in combinatorial optimization problems and serves as a starting point for the possible development new algorithms for RCSPP.

\paragraph{Future research} This work opens new avenues for further investigation. 

First, extending the analysis to a broader range of optimization problems and additional RCSPP instances, including those with more complex and multi-resource profiles, would be valuable. This would naturally require designing new features and metrics to capture the added complexity.

Second, it would be interesting to assess the impact of more complex features that cannot be computed in constant time on prediction accuracy. While a trade-off with computational efficiency is expected, leveraging dynamic information, potentially obtained through pre-processing or heuristic methods, could enhance classification performance.

Lastly, it would be valuable to explore the integration of such predictors within the solution process to potentially develop both new heuristic and exact labeling algorithms. This could lead to more efficient strategies for pruning and guiding the search, ultimately improving solution quality and computational performance.

%\section*{Declarations} 

%The authors have no competing interests to declare that are relevant to the content of this article.

%\paragraph{Data availability} The problem instances used in this study are available online in the SPPRCLIB library. The labels used for training and evaluation were generated using the open-source PathWyse library with the settings described in this article. The complete dataset, including all labels, will be published on Zenodo and made available via DOI upon acceptance of the article.

\newpage
\bibliography{bibliography.bib}
\newpage
\appendix
\appendix
\section{Appendix}
\subsection{Instances of the Dataset} \label{appendix:instances}

Class A
\begin{itemize}
    \item A-n54-k7, A-n60-k9, A-n61-k9, A-n62-k8, A-n63-k9, A-n63-k10,A-n64-k9, A-n65-k9, A-n69-k9, A-n80-k10

\end{itemize}
Class B
\begin{itemize}
    \item B-n45-k6, B-n50-k8, B-n52-k7, B-n66-k9, B-n67-k10, B-n68-k9, B-n78-k10
 
\end{itemize}
Class E
\begin{itemize}
    \item  E-n76-k7, E-n76-k10, E-n76-k14, E-n76-k15, E-n101-k8, E-n101-k14

\end{itemize}
Class M
\begin{itemize}
    \item M-n101-k10, M-n151-k12, M-n200-k17

\end{itemize}
Class P
\begin{itemize}
    \item  P-n50-k7, P-n50-k8, P-n50-k10, P-n51-k10, P-n55-k10, P-n55-k15, P-n55-k7, P-n55-k8, P-n60-k10, P-n60-k15, P-n65-k10, P-n70-k10, P-n76-k4, P-n76-k5, P-n101-k4

\end{itemize}
\subsection{Features of the Dataset} \label{appendix:features}

For each label generated by the optimization algorithm, we collected the following characteristics. We note that all information is computed in constant time, that is, during the execution of the algorithm, we do not perform additional computations to access data, whenever a label is created. In the following, we mark with an asterisk the features that are used to train supervised models in Section \ref{sec:learning}.

\subsubsection*{Algorithm Features}
\begin{itemize}
    \item executionID
    \item iteration
\end{itemize}

For each label, we keep track the ID of the test (executionID) and the iteration of the DSSR-R relaxation that generated it (iteration).

\subsubsection*{General features}

\begin{multicols}{2}
\begin{itemize}
    \item direction
    \item node
    \item predecessor
    \item objective*    
    \item consumption\_critical*
    \item efficiency*
    \item dominated
\end{itemize}
\end{multicols}

We also record if the label was generated with a forward or backward search (direction), the currently reached node (node), its predecessor (precedessor), the total cumulative objective and consumption of critical resource of the partial path (respectively, objective and consumption\_critical). A custom efficiency metric (efficiency), that we defined in section \ref{normalization}.
Finally, we also record if a label is dominated or not. That is, dominated labels are marked with 1 and Pareto labels are marked with 0.

\subsubsection*{Tour-related features}
\begin{multicols}{2}

\begin{itemize}
    \item tour\_length*
    \item nvisited*
    \item nvisited\_unaltered*
    \item nunreachable*
    \item nunreachable\_unaltered*
    \item repeated\_visits*
    \item visited
    \item visited\_unaltered
    \item unreachable
\end{itemize}
\end{multicols}

We then record features related to the nodes visited by the partial path. These include the tour length, the number of loops in the path (repeated\_visits), and both the number of visited (nvisited) and unreachable (nunreachable) nodes, as determined by the current DSSR-R sets. In addition, we store unmodified versions of these features, that is, the actual counts of visited and unreachable nodes that would result without any relaxation.
Finally, we initially save the full bitsets representing the visited and unreachable nodes.

\subsubsection*{Network Status Features}
\begin{multicols}{2}

\begin{itemize}
    \item nlabels\_network*
    \item nlabels\_dominated\_network*
    \item nlabels\_closed\_network*
    \item nlabels\_open\_network*
    \item nlabels\_node*
    \item nlabels\_closed\_node*
    \item nlabels\_open\_node*
\end{itemize}
\end{multicols}

The final set of features we monitor captures the state of the network at the moment each label is generated by the algorithm. Specifically, we record the total number of labels in the network (nlabels\_network), distinguishing between open, closed, and dominated labels. Additionally, for each label, we track the total number of labels in the corresponding node data pool (nlabels\_node), separately identifying open and closed labels. We remark that these features are specific for the direction of the search (i.e. forward or backward).

\subsection{Validation Results} \label{appendix:validation}

We report in Table \ref{tab:validation} the validation results for the classifiers discussed in Section \ref{sec:learning}. For both offline and online tests and each dataset we report the overall accuracy, precision, recall and F1-score. For offline experiments, we also report the training strategy that was adopted.
    
\begin{table}[h!]
\centering
\begin{tabularx}{0.9\textwidth}{X*{7}{r}}\toprule
Test & Strategy & Dataset & Accuracy & Precision & Recall & F1 \\
\hline
\\
\multirow{5}{*}{Offline} 
    & \multirow{2}{*}{one-vs-all} 
        & G & 0.92	&1.00	&0.92	&0.96\\
    &   & I & 0.83	&0.79	&0.82	&0.79\\
    \\
    & \multirow{2}{*}{all-but-one} 
        & G & 0.91	&1.00	&0.91	&0.95 \\
    &   & I & 0.79	&0.88	&0.76	&0.82 \\
    \\
\hline
\\
\multirow{2}{*}{Online} 
    & \multirow{2}{*}{} 
        & G & 0.92	&1.00	&0.92	&0.96 \\
    &   & I & 0.83	&0.79	&0.82	&0.80 \\
\\

\bottomrule
\end{tabularx}

\caption[validation]
{Validation results for every experiment, for models trained on Dataset G and I.} \label{tab:validation}
\end{table}

\subsection{Alternative aggregation strategies for offline learning} \label{appendix:aggregation}
Several alternative aggregation strategies were explored during preliminary experiments.

First, we tested a random-split approach, where one third of the instances were selected at random for training, and the remaining ones were used for testing. This process was repeated multiple times. While this strategy generally outperformed the one-vs-all setup, it underperformed compared to the all-vs-one approach. Moreover, due to the randomized selection, results were highly inconsistent.

Next, we explored class-based aggregation. Since instances within the same class share structural characteristics, such as node distributions and value patterns, we tested whether a model trained on a subset of instances from a given class could effectively predict other instances from the same class. However, results showed that such models were, on average, quite inaccurate.

We also experimented with cross-class predictions: using a model trained on one class to predict instances from other classes. Interestingly, this approach led to better average performance across both datasets compared to the one-vs-all strategy. Nonetheless, some models (notably those trained on classes E and M) displayed a strong bias toward one class. Indeed, these two models require further tuning: better results might be obtained calibrating the Pareto and Dominated labels split used during training, for each instance.

In summary, the models were unable to generalize within the same class. However, as demonstrated by the all-but-one strategy, training on a diverse set of instances, including those from different classes, improves predictive performance.

\end{document}